\documentclass[11pt]{article}
\usepackage{amsmath,amsfonts,amssymb,amsthm,mathrsfs}
\usepackage[a4paper,vmargin={3.5cm,3.5cm},hmargin={2.5cm,2.5cm}]{geometry}
\usepackage[font=sf, labelfont={sf,bf}, margin=1cm]{caption}
\usepackage{graphicx,graphics}
\usepackage{epsfig}
\usepackage{latexsym}
\usepackage[applemac]{inputenc}
\usepackage{ae,aecompl}
\usepackage[english]{babel}
 \usepackage[colorlinks=true]{hyperref}
\usepackage{pstricks}
\usepackage{enumerate}
\usepackage{cleveref}
  \crefname{theorem}{Theorem}{Theorems}
  \crefname{lemma}{Lemma}{Lemmas}
  \crefname{remark}{Remark}{Remarks}
  \crefname{proposition}{Proposition}{Propositions}
  \crefname{definition}{Definition}{Definitions}
  \crefname{corollary}{Corollary}{Corollaries}
  \crefname{section}{Section}{Sections}
  \crefname{figure}{Figure}{Figures}

\newtheorem{theorem}{Theorem}[]
\newtheorem{conjecture}{Conjecture}[]

\newtheorem{proposition}[theorem]{Proposition}
\newtheorem{lemma}[theorem]{Lemma}

\theoremstyle{definition}
\newtheorem*{remark}{Remark}

\title{Planar stochastic hyperbolic infinite triangulations}

\author{Nicolas Curien\thanks{ CNRS and Université Paris 6, E-mail: nicolas.curien@gmail.com} }

\date{}

\begin{document}
\maketitle
\begin{abstract} Pursuing the approach of \cite{AR13} we introduce and study a family of random infinite triangulations of the full-plane that satisfy a natural spatial Markov property. These new random lattices  naturally generalize Angel \& Schramm's Uniform Infinite Planar Triangulation (UIPT) and are hyperbolic in flavor. We prove  that they exhibit a sharp exponential volume growth, are non-Liouville, and that the simple random walk on them has positive speed almost surely.  We conjecture that these infinite triangulations are the local limits of uniform triangulations whose genus is proportional to the size. 
\end{abstract}

\begin{figure}[!h]
 \begin{center}
 \includegraphics[width=9cm]{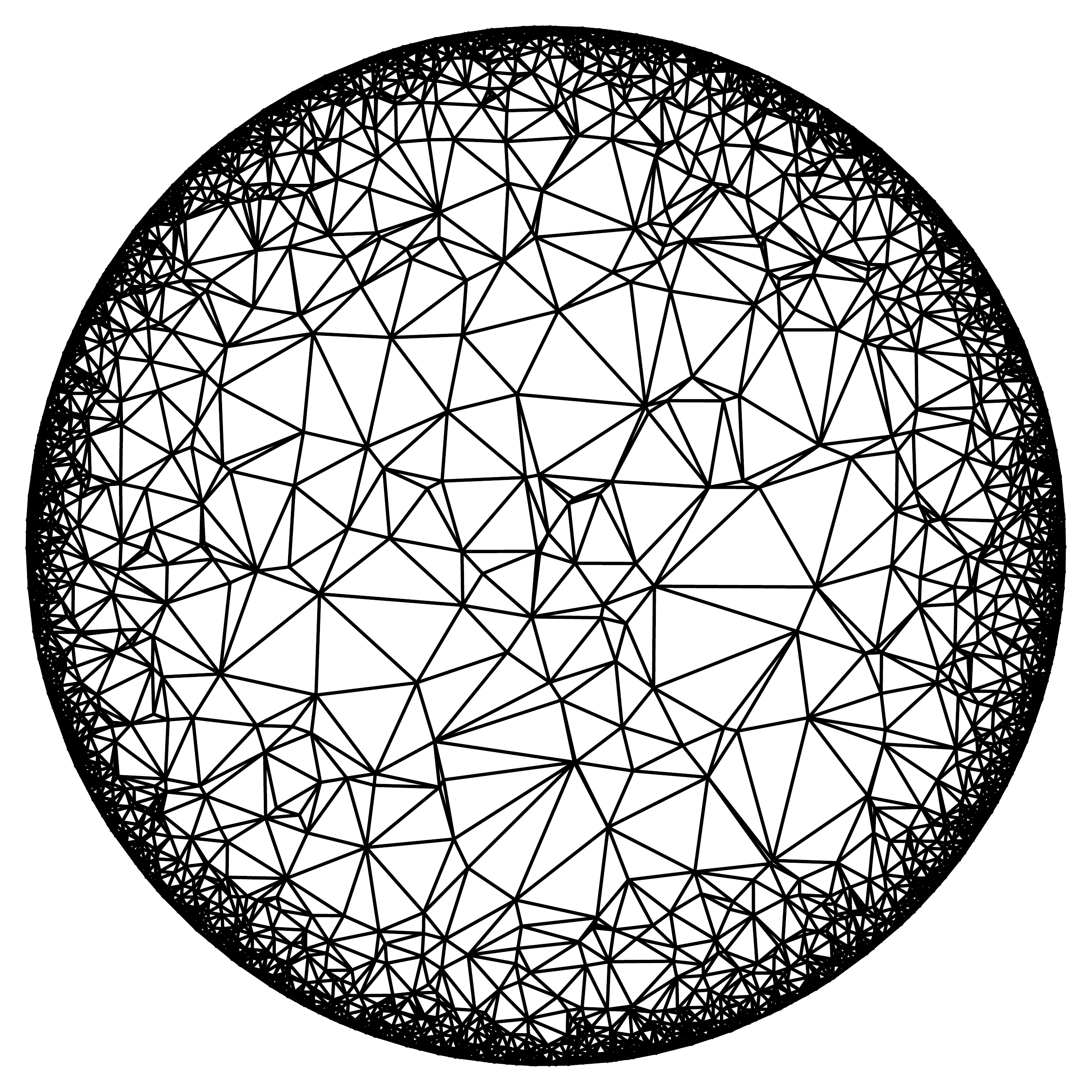}\\ 
 \textsf{An artistic representation of a random ($3$-connected) triangulation of the plane with hyperbolic flavor.}
 \end{center}
 \end{figure}
\medskip

\clearpage

\section*{Introduction} 
Since the introduction of the Uniform Infinite Planar Triangulation by Angel \& Schramm \cite{AS03} as the local limit of large uniform triangulations, a large body of work has been devoted to the study of local limits of random maps and especially random triangulations and quadrangulations, see e.g.\,\,\cite{AS04,ACCR13,BS01,BS13,CD06,Kri05} and the references therein. Recently, Angel \& Ray \cite{AR13} classified all random triangulations of the half-plane that satisfy a natural spatial Markov property and discovered new random lattices of the half-plane exhibiting  a ``hyperbolic'' behavior \cite{Ray13}. Motivated by these works we construct the analogs of these lattices in the full-plane topology and study their properties in fine details. 

\paragraph{Classification of Markovian triangulations of the plane.} Recall that a planar map is a proper embedding of a finite connected graph into the sphere seen up to deformations that preserve the orientation.  All maps considered here are rooted, i.e.~given with a distinguished oriented edge. A triangulation is a planar map whose faces have all degree three. For the sake of simplicity we restrict ourselves to $2$-connected triangulations where loops are forbidden (but multiple edges are allowed). We shall also deal with \emph{infinite triangulations of the plane} or equivalently infinite triangulations with one end (a graph $G$ is said to have one end if $G \backslash H$ contains exactly one infinite connected component for any subgraph $H \subset G$). They can be realized as proper embeddings  $ \mathcal{X}$ (seen up to continuous deformations preserving the orientation)   of graphs in $ \mathbb{R}^2$ such that every compact subset of $ \mathbb{R}^2$ intersects only finitely many edges of $ \mathcal{X}$ and such that the faces have all degree three. For any $p \geq 2$, a triangulation $t$ of the $p$-gon, also called triangulation with a boundary of perimeter $p$, is  a (rooted) planar map whose faces are all triangles except for one distinguished face of degree $p$, called the hole, whose boundary is made of a simple cycle (no pinch point). The perimeter $|\partial t|$ of $t$ is the degree of its hole and its size $|t|$ is its number of  vertices. The set of all finite triangulations of the $p$-gon is denoted by $ \mathcal{T}_{p}$ and we set $ \mathcal{T}_{B} = \cup_{p \geq 2} \mathcal{T}_{p}$ for the set of all finite triangulations with a boundary. \medskip

 Here comes the key definition of this work: We say that a random infinite triangulation $ \mathbf{T}$ of the plane is \emph{$\kappa$-Markovian} for $\kappa >0$ if there exist  non-negative numbers $(C_{i}^{(\kappa)} : i \geq 2)$ such that for any $t \in \mathcal{T}_{p}$  we have   \begin{eqnarray} P( t \subset  \mathbf{T}) &=& C_{p} ^{(\kappa)} \cdot \kappa^{| t|},  \label{eq:SMP}\end{eqnarray} where by $ t \subset  \mathbf{T}$ we mean that $ \mathbf{T}$ is obtained from $t$ (with coinciding roots) by filling its hole with a necessary unique infinite triangulation of the $p$-gon.  Our first result which parallels \cite{AR13} is to show the existence and uniqueness of a one-parameter family of such triangulations:
\begin{theorem} \label{thm:un}For any $\kappa \in (0, \frac{2}{27}]$ there exists a unique (law of a) random $\kappa$-Markovian triangulation $ \mathbf{T}_{\kappa}$ of the plane. If $ \kappa >\frac{2}{27}$ there is none.
\end{theorem}

 In the special case $\kappa = \frac{2}{27}$, called the critical case, it follows from \cite[Theorem 5.1]{AS03} that the triangulation $ \mathbf{T}_{{2}/{27}}$ has the law of the uniform infinite planar triangulation (UIPT) introduced by Angel \& Schramm \cite{AS03} as the limit of uniform triangulations of the sphere of growing sizes. The UIPT and its quadrangular analog the UIPQ have received a lot of attention in recent years \cite{Ang03,BCsubdiffusive,CMMinfini,GGN12,LGM10} partially motivated by the connections with the physics theory of 2-dimensional quantum gravity, the Gaussian free field \cite{CurKPZ} and the Brownian map \cite{CLGplane,LG11,Mie11}. Many fundamental problems about the UIPT/Q are still open. We will see below that the qualitative behavior of $ \mathbf{T}_{\kappa}$ in the regime $\kappa < \frac{2}{27}$, called hyperbolic regime (this terminology will be justified by the following results), is much different  from that of the UIPT. 
 \medskip 
 
   In the work  \cite{AR13}, the authors classified all random triangulations of the \emph{half-plane} that satisfy a very natural, but slightly different, spatial Markov property: a random triangulation of the half-plane has the spatial Markov property of \cite{AR13} if conditionally on any simply connected neighborhood of the root (necessarily located on the infinite boundary), the remaining lattice has the same law as the original one. Angel \& Ray classified these lattices using a single parameter $\alpha \in [0,1)$ which is equal to the probability that the face adjacent to a given edge on the boundary is a triangle pointing inside the map. Our lattices $ (\mathbf{T}_{\kappa})$ for $\kappa \in (0,\frac{2}{27}]$   are the full-plane analogs of the half-planar lattices of \cite{AR13} for $\alpha \geq 2/3$ and 
 \begin{eqnarray} \frac{\alpha^2(1-\alpha)}{2} =\kappa \qquad \mbox{and} \qquad \alpha \in [2/3,1).   \label{eq:alphakappa}\end{eqnarray}
Angel \& Ray also exhibited subcritical half-planar lattices corresponding to $\alpha \in (0,2/3)$ which are tree-like \cite{Ray13}. In our full-plane setup, no such subcritical phase exists. Although similar in spirit to Angel \& Ray's spatial Markov property our Markovian assumption \eqref{eq:SMP} is slightly different mainly because of the topology of the plane which forces the presence of a function of the perimeter $(C_{p}^{(\kappa)}: p \geq 2)$ and also because we impose an exponential dependence in the size.  

\paragraph{Peeling process.} In the theory of random planar maps, the spatial Markov property is a key feature that has already been thoroughly used, generally under the form of the ``peeling process'' \cite{Ang03,ACpercopeel,BCsubdiffusive,CurKPZ,CLGpeeling,MN13}. The peeling process has been conceived by Watabiki \cite{Wat95} and formalized by Angel \cite{Ang03} in the case of the UIPT. This is an algorithmic procedure that enables to construct the lattice in a Markovian fashion by exploring it face after face (possibly revealing the finite regions enclosed). It turns out that equation \eqref{eq:SMP} implies that $ \mathbf{T}_{\kappa}$ must admit such a peeling process which yields a path to prove both existence and uniqueness in Theorem \ref{thm:un}, as in \cite{AR13}. Furthermore,  we will see in Section \ref{sec:harmonic} that the function $p \mapsto C_{p}^{(\kappa)}$ of \eqref{eq:SMP} will be interpreted as a harmonic function of the underlying random walk governing the construction of $ \mathbf{T}_{\kappa}$ by the peeling process. The peeling process is also a key tool in the proof of the up-coming Theorems \ref{thm:volume} and \ref{thm:hyperbolic}.\medskip

\paragraph{Properties of the planar stochastic hyperbolic infinite triangulations.}
Let us now turn to the properties of these new random lattices in the hyperbolic regime $ \kappa \in (0 , \frac{2}{27})$.  If $ \mathsf{T}$ is a finite triangulation or an infinite triangulation of the plane, we let $B_{r}( \mathsf{T})$ denote the subtriangulation obtained by keeping the faces of $ \mathsf{T}$ that contain at least one vertex at graph distance less than or equal to $r-1$ from the origin of the root edge in $ \mathsf{T}$. Hence, $B_{r}( \mathsf{T})$ is a triangulation with a finite number of holes. In the infinite case, we also denote by $ \overline{B}_{r}( \mathsf{T})$ the hull of the ball obtained by filling-in all the finite components of $  \mathsf{T} \backslash B_{r}( \mathsf{T})$. Since $ \mathsf{T}$ is one-ended, $\overline{B}_{r}( \mathsf{T})$ belongs to $ \mathcal{T}_{B}$ and its boundary in $ \mathsf{T}$ is a simple cycle made of edges whose vertices are at distance exactly $r$ from the origin of the root edge in $ \mathsf{T}$.

\begin{theorem}[Sharp exponential volume growth] \label{thm:volume} For any $\kappa \in (0, \frac{2}{27})$ introduce $\alpha \in ( \frac{2}{3},1)$ satisfying \eqref{eq:alphakappa}  and let $ \delta_{\kappa} = \sqrt{\alpha(3\alpha - 2)}$.  There exists a random variable $ \Pi_{\kappa}$ such that $\Pi_{\kappa} \in (0, \infty)$ almost surely with
 \begin{eqnarray*} \left( \frac{\alpha - \delta_{\kappa}}{\alpha+ \delta_{\kappa}}\right)^n |\partial \overline{B}_{n}( \mathbf{T}_{\kappa})| & \xrightarrow[n\to\infty]{a.s.}&  \Pi_{\kappa},\\
\mbox{ and} \qquad  \frac{ | \overline{B}_{n}(\mathbf{T}_{\kappa})|}{| \partial \overline{B}_{n}( \mathbf{T}_{\kappa})|} &\xrightarrow[n\to\infty]{a.s.}& \frac{\alpha(2 \alpha -1)}{\delta_{\kappa}^2} .  \end{eqnarray*}
\end{theorem}
In \cite{Ray13}, exponential bounds for the volume growth in the half-plane version of $ \mathbf{T}_{\kappa}$ are obtained but with non-matching exponential factors. It is very likely that the methods used here can be employed to settle \cite[Question 6.1]{Ray13}. The results of the last theorem should be compared with the analogous properties for supercritical Galton--Watson trees (whose offspring distribution satisfies the $x\log x$ condition) where the number of individuals at generation $n$, properly normalized, converges towards a non-degenerate random variable on the event of non-extinction. We also show that in the hyperbolic regime, $ \mathbf{T}_{\kappa}$ has a positive anchored expansion constant (Proposition \ref{prop:anchored}). \medskip

We then turn to the study of the simple random walk on $ \mathbf{T}_{\kappa}$: Conditionally on $\mathbf{T}_{\kappa}$ we launch a random walker from the target of the root edge and let it choose inductively one of its adjacent oriented edges for the next step. In the critical case $\kappa = \frac{2}{27}$, the simple random walk on the UIPT is known to be recurrent \cite{GGN12}. We show here that the behavior of the simple random walk is drastically different when $\kappa < \frac{2}{27}$. Recall that a graph is non-Liouville if and only if it possesses non-constant bounded harmonic functions.

\begin{theorem}[Hyperbolicity] \label{thm:hyperbolic} For $\kappa \in (0, \frac{2}{27})$ there exists $s_{\kappa} >0$ such that almost surely
$$\lim_{n \to \infty} n^{-1} \mathrm{d_{gr}}(X_{0},X_{n}) = s_{\kappa},$$ where $(X_{i})_{i\geq 0}$ are the vertices visited by the simple random walk and $ \mathrm{d_{gr}}$ is the graph metric. Also, $ \mathbf{T}_{\kappa}$ is almost surely non-Liouville in the hyperbolic regime.
\end{theorem}

The connoisseurs may remember that a major difficulty towards proving the recurrence of the UIPT \cite{GGN12} was the lack of a uniform bound on the degree. The situation is similar here, as positive speed would directly follow from positive anchored expansion in the bounded degree case by the result of \cite{Vir00}. The unboundedness of the vertex degrees in $ \mathbf{T}_{\kappa}$ forces us to find a different technique. The proof of Theorem \ref{thm:hyperbolic} occupies the major part of Section \ref{sec:speed} and makes extensive use of the fact that the random lattice $ \mathbf{T}_{\kappa}$ is stationary and reversible with respect to the simple random walk (Proposition \ref{prop:reversible}). In words, re-rooting $ \mathbf{T}_{\kappa}$ along a simple random walk path does not change its distribution. This stochastic invariance by translation replaces the deterministic invariance of transitive lattices. This is a key feature of the full-plane models compared to the half-plane models of \cite{AR13}.  The proof of Theorem \ref{thm:hyperbolic} also combines several geometric arguments such as: the exploration process of $ \mathbf{T}_{\kappa}$ along the simple random walk of \cite{BCsubdiffusive}, the recent results of \cite{BCGLiouville} on intersection properties of planar lattices and the entropy method for stationary random graphs \cite{BCstationary}. 

\paragraph{A speculation.} We end this introduction by stating a conjecture relating our planar stochastic hyperbolic infinite triangulations to local limits of triangulations in high genus. More precisely, let $ \mathcal{T}_{n,g} $ be the set of all (rooted) triangulations of the torus of genus $g \geq 0$ with $n$ vertices and denote by $T_{n,g}$ a random uniform element in $ \mathcal{T}_{n,g}$. We conjecture that there is a continuous decreasing function function $f(\theta) \in (0, \frac{2}{27}]$ with $f(\theta) \to 0$ as $\theta \to \infty$ and $f(0) = \frac{2}{27}$ such that we have the following convergence
$$ T_{n,[\theta n]} \quad \xrightarrow[n\to\infty]{(d)} \quad \mathbf{T}_{f(\theta)},$$ for the local topology. 
See Section \ref{sec:comments} for a more precise statement. Some results in the literature already indicate that large triangulations of high genus should be locally planar such as the work of Guth, Parlier and Young on pants decomposition of random surfaces \cite{GPY11} or the paper \cite{ACCR13} on the case of unicellular maps. \medskip

The organization of the paper should be clear from the table of contents below. \bigskip

\noindent \textbf{Acknowledgments:} I am grateful to Omer Angel, Itai Benjamini, Guillaume Chapuy and Gourab Ray for useful discussions on and around Conjecture 1.

\tableofcontents

\clearpage

\section{Construction of $( \mathbf{T}_{\kappa})$ for  $\kappa \in (0, \frac{2}{27}]$}
Fix $\kappa >0$. We will first show that the law (if it exists) of a $\kappa$-Markovian random infinite triangulation of the plane is unique and that $\kappa$ must be less than or equal to $\frac{2}{27}$. In  this whole section we thus assume the existence of $ \mathbf{T}_{\kappa}$, a $\kappa$-Markovian triangulation  of the plane.
\subsection{Uniqueness} \label{sec:uniqueness}
We start with a few pre-requisites on the local topology on triangulations. Following \cite{BS01}, if $t$ and $t'$ are two rooted finite triangulations, the local distance between $t$ and $t'$ is set to be
 \begin{eqnarray*} \mathrm{d_{loc}}(t,t') &=& \big(1+ \sup\{r\geq 0 : B_{r}(t)=B_{r}(t')\}\big)^{-1}.\end{eqnarray*}
The set of all finite triangulations is not complete for this distance and we shall add infinite triangulations to it. The metric space $( \mathcal{T}_{\infty}, \mathrm{d_{loc}})$ we obtain is then Polish. Since $ \mathbf{T}_{\kappa}$ is an infinite random triangulation \emph{with only one end}, it is easy to see that its law in $( \mathcal{T}_{\infty}, \mathrm{d_{loc}})$ is characterized by the values of $P(t \subset \mathbf{T}_{p})$ for all triangulations with a boundary $t \in \mathcal{T}_{B}$.  Hence, establishing the uniqueness of the law of $ \mathbf{T}_{\kappa}$ reduces to showing that the function $(C_{i}^{(\kappa)}:i \geq 2)$ involved in \eqref{eq:SMP} is uniquely characterized by $\kappa >0$.

We begin with a  simple remark. Since $ \mathbf{T}_{\kappa}$ is a $2$-connected triangulation, the triangle on the left of the root edge is necessarily a triangle with $3$ distinct vertices with the root edge located on one of its side that we see as a triangulation of the $3$-gon denoted by $t_{0}$. By \eqref{eq:SMP} we must have 
 \begin{eqnarray} \label{eq:C3=1} 1 =P(t_{0}\subset  \mathbf{T}_{ \kappa}) = \kappa^3 C_{3}^{(\kappa)}.  \end{eqnarray}

\begin{figure}[!h]
 \begin{center}
 \includegraphics[width=12cm]{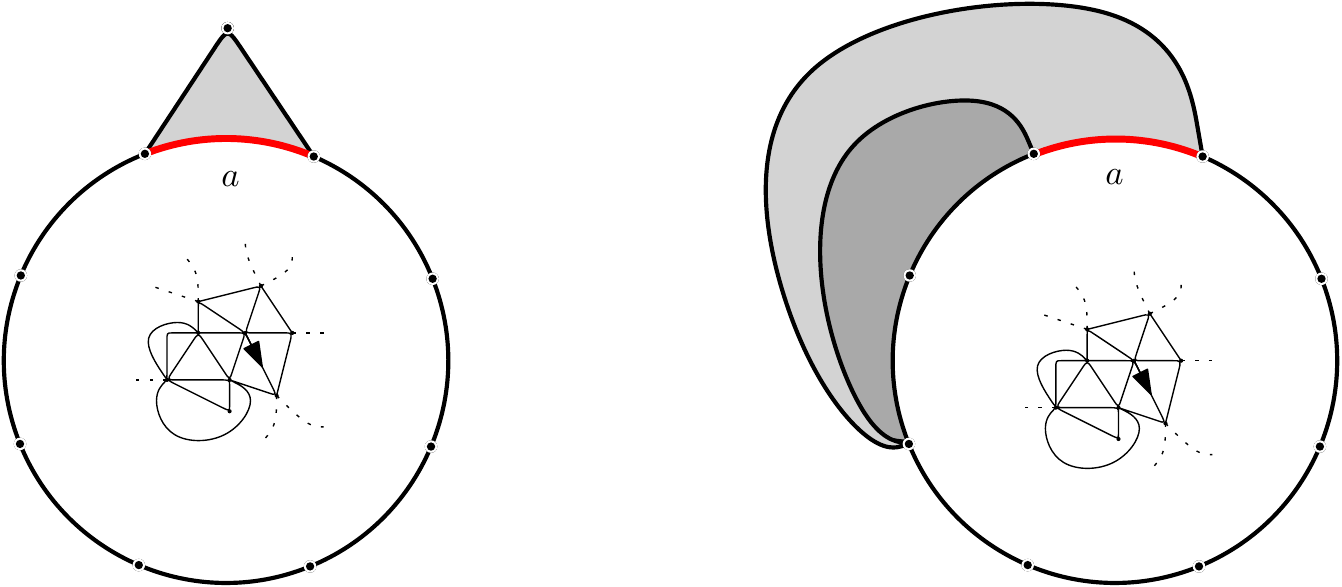}
 \caption{ \label{fig:peeling}Peeling the edge $a$: the triangle revealed is in light gray and the finite enclosed region is in dark gray on the second figure.}
 \end{center}
 \end{figure}

We will now get another relation linking $ \kappa$ and the $(C_{p}^{(\kappa)})_{p \geq 1}$. This is done by increasing a map using the so-called peeling mechanism. Let $ \mathsf{T}$ be a triangulation of the plane and assume that  $t \subset \mathsf{T}$ for some  $t  \in \mathcal{T}_{p}$. For any edge $a$ on the boundary of $t$ we condiser the triangulation which is obtained by adding to $t$ the triangle adjacent to $a$ in $ \mathsf{T} \backslash t$ as well as the finite region this triangle may enclose (recall that $ \mathsf{T}$ is one-ended). We call this operation peeling the edge $a\in \partial t$. Two different situations may appear : either the triangle revealed contains a vertex inside $ \mathsf{T}\backslash t$ (left on Fig.~\ref{fig:peeling}) or this triangle ``swallows'' $k$ edges on the boundary of $t$ either to the left or to the right of $a$ and encloses a finite triangulation of the $k+1$-gon (right on Fig.~\ref{fig:peeling}). Note that $k \in \{1,2, ... , p-2\}$ where $p$ is the perimeter of $t$. 

For two triangulations with a boundary $t$ and $t'$ we write $(t,a) \to t'$ if $t'$ is a possible outcome of the peeling of the edge $a \in   \partial t$ in some underlying triangulation $ \mathsf{T}$. It is easy to see that such $t'$ are  obtained by either gluing triangle to $a$ outside $t$ or by gluing a triangle to $a$ with its third vertex identified with a vertex of the boundary of $t$ and filling one of the two holes created with a finite triangulation having the proper perimeter. In the first case the size of the triangulation increases by $1$, and in the second case it increases by the number of inner vertices (not located on the boundary) of the enclosed triangulation. This operation is rigid in the sense that two different ways of increasing $t$ yield to two distinct maps. For $p \geq 2$ and $n \geq p$, pick  a triangulation $t$ of the $p$-gon with $n$  vertices and fix deterministically an edge $a$ on its boundary. By the previous discussion we have  
 \begin{eqnarray} C_{p}^{(\kappa)} \cdot \kappa^{n} \quad \underset{ \eqref{eq:SMP}}{=} \quad  P( t \subset \mathbf{T}_{ \kappa}) &=& P \left(\bigcup_{ (t,a) \to t'} \{ t' \subset \mathbf{T}_{\kappa}\}\right)\nonumber \\ &\underset{ \mathrm{rigidity}}{=}& \sum_{ (t,a) \to t'} P\left( t' \subset \mathbf{T}_{\kappa}\right)\nonumber \\	
 &\underset{\eqref{eq:SMP}}{=}& C_{p+1}^{(\kappa)}\kappa^{n+1} + 2\sum_{i=1}^{p-2} C_{p-i}^{(\kappa)}  \cdot \kappa^{n} \sum_{\tau \in \mathcal{T}_{i+1}} \kappa^{|{\tau}|-i-1} \label{eq:interm}.
  \end{eqnarray} 
Now, for $i \geq 1$ and $\kappa >0$ introduce the functions $$Z_{i+1}^{(\kappa)} = \displaystyle \sum_{\tau \in \mathcal{T}_{i+1}}  \kappa^{|\tau|-i-1} = \sum_{\tau \in \mathcal{T}_{i+1}}  \kappa^{|\tau|-| \partial \tau|}.$$ A closed formula is known for these numbers (see \cite{GJ83} or \cite[Proposition 2.4]{AS03}) and they are finite if and only if $\kappa \in (0 , \frac{2}{27}]$.  They can be interpreted as the partition function of the following probability measure: The \emph{Boltzmann} probability distribution of the $i+1$-gon with parameter $\kappa$ is the probability measure that assigns weight $\kappa^{|t|-i-1}/Z_{i+1}^{(\kappa)}$ to each triangulation $t$ of the $i+1$-gon. Hence \eqref{eq:interm} becomes
    \begin{eqnarray} \label{eq:cplinked}
    C_{p}^{(\kappa)}  = \kappa \cdot C_{p+1}^{(\kappa)} + 2 \sum_{i=1}^{p-2}  C_{p-i}^{(\kappa)} \cdot Z_{i+1}^{(\kappa)} \qquad \forall p \in \{2,3,4,...\}.  \end{eqnarray}
If $\kappa >\frac{2}{27}$  then $Z_{i}^{(\kappa)}=\infty$ for any $i \geq 2$, so we must suppose that $\kappa \in (0, \frac{2}{27}]$. Using the last display with $p=2$ we find that $C_{2}^{(\kappa)} = \kappa C_{3}^{(\kappa)}$ which combined with \eqref{eq:C3=1} fixes the value of $C_{2}^{(\kappa)}$. Next, using \eqref{eq:cplinked} recursively for $p=3, 4, \ldots$ we see that the values of $C_{p}^{(\kappa)}$ for $p \geq 4$ are fixed by $\kappa$ only. This proves \emph{uniqueness} of the law of $ \mathbf{T}_{\kappa}$.

\subsection{Interpreting  $(C_{p}^{(\kappa)} : p \geq 2)$ as a harmonic function}
\label{sec:harmonic}
Fix $\kappa \in (0, \frac{2}{27}]$ and \emph{define} the numbers $C_{2}^{(\kappa)}, C_{3}^{(\kappa)},...$ using \eqref{eq:C3=1} and \eqref{eq:cplinked} as in the preceding section. Towards proving the existence of a $\kappa$-Markovian triangulation, our first task is to show that   \begin{eqnarray} C_{p}^{(\kappa)} > 0\quad \mbox{ for every } p \geq 2.  \end{eqnarray} To do so it will be very useful to interpret them probabilistically as in \cite{CLGpeeling}. We start by recalling a key calculation that can be found in  \cite[Section 3.1]{AR13}. Let $\alpha \in [2/3,1)$ given by \eqref{eq:alphakappa} and let $\beta = \kappa/\alpha$, then we have $$1 = \alpha + 2 \sum_{i=1}^\infty \beta^i Z_{i+1}^{(\alpha \beta)}.$$
This enables us to define a probability distribution $\boldsymbol{q}^{(\kappa)} = \{ ... , q_{-3}^{(\kappa)},q_{-2}^{(\kappa)},q_{-1}^{(\kappa)}, q_{1}^{(\kappa)}\} $ by setting $$q^{(\kappa)}_{1}  = \alpha \quad  \mbox{and} \quad   q^{(\kappa)}_{-i} = 2 \beta^i Z_{i+1}^{(\kappa)} \quad \mbox{ for }i \geq 1.$$
From \cite[Equation (3.7)]{AR13} we even have an exact formula 
 \begin{eqnarray} \label{eq:exact}q^{(\kappa)}_{-i} = \frac{2}{4^i}\frac{(2i - 2)!}{(i - 1)! (i + 1)!}\left(\frac{2}{\alpha} - 2\right)^i\big((3\alpha - 2)i + 1\big), \quad \mbox{for } i \geq 1.  \end{eqnarray}
Finally we introduce $(\Xi_{n}^{(\kappa)})_{n \geq 0}$  a random walk started from $2$ with independent increments following the distribution $ \boldsymbol{q}^{(\kappa)}$. A computation using \eqref{eq:exact} done in \cite[Lemma 4.2]{Ray13} shows that the drift of this walk is   given by
  \begin{eqnarray}  \delta_{\kappa}:= \sum_{i \leq 1} i q_{i}^{(\kappa)} =  \sqrt{\alpha(3\alpha - 2)} >0. \label{eq:mean} \end{eqnarray}
We can now show that $C_{p}^{(\kappa)} >0$ for all $p \geq 2$: In \eqref{eq:cplinked} we multiply both sides by $\beta^p$ and set $\tilde{C}_{p}^{(\kappa)} = \beta^p C_{p}^{(\kappa)}$ for $p \geq 2$ and put $ \tilde{C}_{p}^{(\kappa)} =0$ otherwise, so that \eqref{eq:cplinked} becomes
 \begin{eqnarray} \label{eq:htransform} \tilde{C}_{p}^{(\kappa)} = \sum_{i \in \{...,-3,-2,-1,1\}}q_{i}^{(\kappa)} \cdot  \tilde{C}_{p+i}^{(\kappa)} \qquad \mbox{ for } p\geq 2.  \end{eqnarray}
In other words, the function $p \mapsto \tilde{C}_{p}^{(\kappa)}$ is (the only) function which is harmonic for the random walk $\Xi^{(\kappa)}$ on $\{2,3,...\}$ and null for $p \leq 1$ subject to the condition $\tilde{C}_{3}^{(\kappa)} = \alpha^{-3}$ given by \eqref{eq:C3=1}. Note that we have $ \tilde{C}_{2}^{(\kappa)} = \alpha  \tilde{C}_{3}^{(\kappa)} < \tilde{C}_{3}^{(\kappa)}$ and that  the last display can be written as $$ q^{(\kappa)}_{1} \big(\tilde{C}_{p+1}^{(\kappa)}-\tilde{C}_{p}^{(\kappa)}\big) = \sum_{i=1}^\infty q^{(\kappa)}_{-i} \big( \tilde{C}_{p}^{(\kappa)} - \tilde{C}_{p-i}^{(\kappa)}\big).$$ We immediately conclude by induction on $p \geq 2$ that  $\tilde{C}_{p}^{(\kappa)}$ is increasing in $p$ and so $ \tilde{C}_{p}^{(\kappa)}$ is positive for all $p$'s. It follows that $C_{p}^{(\kappa)} >0$ for every $p \geq 2$ as desired. 

In the case $\kappa= \frac{2}{27}$ (equivalently $\alpha= \frac{2}{3}$) the functions $ \tilde{C}_{p}^{(2/27)}$ are explicitly known and correspond to the function $9^{-p}C_{p}$ in  \cite{Ang03} and thus grows like $\sqrt{p}$ when $p \to \infty$. In the hyperbolic regime a different behavior appears:

\begin{lemma} \label{lem:fini} When $\kappa < \frac{2}{27}$, the increasing sequence $ \tilde{C}_{p}^{(\kappa)}$ converges to $ (\alpha \delta_{\kappa})^{-1}$ as $p \to \infty$.
\end{lemma}
\proof By monotonicity $\lim_{p \to \infty} \tilde{C}_{p}^{(\kappa)}$ exists in $(0, \infty]$. Let $\Xi^{(\kappa)}$ be the random walk started from $2$ with i.i.d.\,\,increments distributed as $ \boldsymbol{q}^{(\kappa)}$. Since its drift $\delta_{\kappa}$ is positive, the stopping time $\tau_{2} = \inf\{ i \geq 0 : \Xi_{i}^{(\kappa)} < 2\}$ has a positive chance to be infinite. For $\tilde{C}_{p}^{(\kappa)}$ is harmonic on $\{2,3,4, ...\}$, the process $\tilde{C}_{\Xi^{(\kappa)}_{n\wedge \tau_{2}}}^{(\kappa)}$ is a martingale and thus 
 \begin{eqnarray*} \tilde{C}_{2}^{(\kappa)} \quad = \quad  E[\tilde{C}_{\Xi^{(\kappa)}_{n}}^{(\kappa)} \mathbf{1}_{\tau_{2} > n}]  \quad  \xrightarrow[n\to\infty]{} \quad  \lim_{p \to \infty} \tilde{C}_{p}^{(\kappa)} P(\tau_{2} = \infty).  \end{eqnarray*}To finish the proof and compute $P( \tau_{2} =\infty)$ we remark that the random walk $\Xi^{(\kappa)}$ has increments bounded above by $1$ so that we can apply the ballot theorem \cite[Theorem 2]{ABR08}: if $\xi_{0}^{(\kappa)}$ are i.i.d.\,\,copies of law $ \boldsymbol{q}^{(\kappa)}$ we have
 \begin{eqnarray*} P( \tau_{2}= \infty) &=& P(2 + \xi_{1}^{(\kappa)}+ \xi_{2}^{(\kappa)} + ... +\xi_{i}^{(\kappa)} \geq 2, \forall i \geq 1)\\
 &=& P(1+ \xi_{0}^{(\kappa)} + \xi_{1}^{(\kappa)}+ \xi_{2}^{(\kappa)} + ... +\xi_{i}^{(\kappa)} \geq 2, \forall i \geq 1 \mid \xi_{0}^{(\kappa)}=1)\\
 &=& \frac{P( \xi_{0}^{(\kappa)}+ \xi_{1}^{(\kappa)} + ... + \xi_{i}^{(\kappa)} > 0, \forall i \geq 1)}{P(\xi_{0}^{(\kappa)}=1)} 
 = \frac{\delta_{\kappa}}{\alpha}.  \end{eqnarray*}\endproof 


\subsection{Peeling construction}

\label{sec:construction}
We now construct the desired lattices $ \mathbf{T}_{\kappa}$. The method is mimicked from  \cite{AR13} and the idea is to revert the procedure used in Section \ref{sec:uniqueness} in order to provide an algorithmic device called the \emph{peeling process} \cite{Ang03} that constructs a sequence of growing triangulations with a boundary.  For a particular peeling procedure, these triangulations are shown to exhaust the plane and define an infinite triangulation with one end.\medskip 

\paragraph{General peeling.} The peeling process depends on an algorithm $ \mathcal{A}$ which associates with every triangulation $t \in \mathcal{T}_{B}$ one of its boundary edges. From this, we construct a growing sequence of triangulations with a boundary $( \mathrm{T}_{n}^{(\kappa), \mathcal{A}} : n \geq 0)$ as follows. To start with, $ \mathrm{T}_{0}^{(\kappa), \mathcal{A}}$ is the root triangulation composed by a single oriented edge (seen as a triangulation of the $2$-gon). Inductively, assume that  $\mathrm{T}_{n}^{(\kappa), \mathcal{A}}$ is constructed. We write $p = |\partial \mathrm{T}_{n}^{(\kappa), \mathcal{A}}|$ and denote by $a\in \partial \mathrm{T}_{n}^{(\kappa)}$ the edge chosen by the algorithm. Notice that this choice may depend on an other source of randomness. Independently of $ \mathrm{T}_{n}^{(\kappa), \mathcal{A}}$ and of the possible extra randomness of $ \mathcal{A}$, the next triangulation $\mathrm{T}_{n+1}^{(\kappa), \mathcal{A}}$ is obtained as follows: With probability $$ q_{1,p}^{(\kappa)}:=q^{(\kappa)}_{1} \cdot  \frac{\tilde{C}_{p+1}^{(\kappa)}}{\tilde{C}_{p}^{(\kappa)}},$$ the triangulation $\mathrm{T}_{n+1}^{(\kappa), \mathcal{A}}$ is obtained from $\mathrm{T}_{n}^{(\kappa), \mathcal{A}}$ by gluing a triangle onto the edge $a$ as in Fig.~\ref{fig:peeling} left. Otherwise, for $-p+2 \leq i \leq  p-2$ with probability 
$$   \frac{1}{2}q_{-|i|,p}^{(\kappa)}:=\frac{1}{2}q^{(\kappa)}_{-|i|}  \cdot \frac{\tilde{C}_{p-|i|}^{(\kappa)}}{\tilde{C}_{p}^{(\kappa)}},$$ we glue a triangle on $ a$  and identify its third vertex with the $|i|$th vertex on the left or on the right of $a$ depending on the sign of $i$ as in Fig.~\ref{fig:peeling} right. Finally,  independently of these choices fill the hole created by the triangle with an independent Boltzmann triangulation of the $|i|+1$-gon with parameter $\kappa$ to get $ \mathrm{T}_{n+1}^{(\kappa), \mathcal{A}}$. According to \eqref{eq:htransform} these probability transitions sum-up to $1$. \medskip 

Assume now that the algorithm $ \mathcal{A}$ is deterministic, i.e.\,\,can be seen as a function $ \mathcal{A}(t) \in \mathrm{Edges}(\partial t)$. Using the same calculations as in Section \ref{sec:uniqueness} one sees by induction that for every $n\geq 0$ and for every triangulation $t \in \mathcal{T}_{p}$ that is a possible outcome of the construction at step $n$ (that is $ P(\mathrm{T}_{n}^{(\kappa), \mathcal{A}} = t) >0$) we have 
  \begin{eqnarray} \label{eq:fixedtime} P( \mathrm{T}_{n}^{(\kappa), \mathcal{A}} = t) = C_{p}^{(\kappa)} \cdot \kappa^{|t|}.  \end{eqnarray} Remark that the right-hand side of the last display does not depend on the order in which the peeling steps are performed nor on $ \mathcal{A}$ as long as $P(\mathrm{T}_{n}^{(\kappa), \mathcal{A}} = t) >0$. The last display can also be extended to the case when $n$ is replaced by a stopping time $\tau$, that is, a random variable such that $ \{\tau=n\}$ is a measurable function of $ \mathrm{T}_{n}^{(\kappa), \mathcal{A}}$. 
  
\begin{figure}[!h]
 \begin{center}
 \includegraphics[width=16cm]{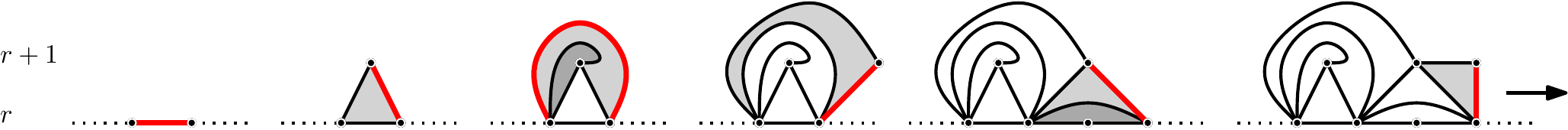}
 \caption{Illustration of the peeling by layers algorithm. \label{fig:layers}}
 \end{center}
 \end{figure}

\paragraph{Defining $ \mathbf{T}_{\kappa}$.}  However, the law of the structure of $( \mathrm{T}_{n}^{(\kappa), \mathcal{A}})$ does depend on the algorithm $ \mathcal{A}$ and it could happen that the increasing union $ \cup_{n \geq 0} \mathrm{T}_{n}^{(\kappa), \mathcal{A}}$ does not create a triangulation of the plane  (imagine for example that one edge is never peeled). To prevent this, we now pick a particular \emph{deterministic} algorithm called $  \mathcal{L}$ for ``layers''. Specifically, peel the left hand side and then the right hand side of the root edge during the first two steps and then, at step $n$, peel the right-most edge on $ \partial \mathrm{T}_{n}^{(\kappa), \mathcal{L}}$ which belongs to the  triangle we just revealed.  See Fig.~\ref{fig:layers}.
 
We easily prove by induction that this algorithm associates with every triangulation $t \in \mathcal{T}_{B}$ an edge $\mathcal{L}(t) \in \partial t$ containing one endpoint $x$ which minimizes $ \left\{ \mathrm{d}_{ \mathrm{gr}}^{ t}(x, e^-) : x \in \partial t \right\},$ where $e^-$ is the origin of the root edge and $ \mathrm{d}_{ \mathrm{gr}}^t(\cdot,\cdot)$ is the graph distance inside $t$. First, an argument similar to \cite[Proposition 3.6]{AR13} or \cite[Lemma 4.4]{Ray13} shows that using this peeling construction, every vertex on the boundary of the growing triangulations will be eventually be swallowed in the process and so 
$$ {\mathbf{T}}_{\kappa} \quad \underset{ \mathrm{def.}}{:=} \quad \bigcup_{n \geq 0} \mathrm{T}_{n}^{(\kappa), \mathcal{L}},$$
   defines an infinite triangulation \emph{of the plane}. We will now check that this random lattice is $\kappa$-Markovian. If $\tau_{r}<\infty$ is the first time when no vertex on $ \partial \mathrm{T}_{n}^{(\kappa), \mathcal{L}}$ is at distance less than  $r-1$ from $e^{-}$  then an easy geometric argument   (see \cite[Proposition 6]{BCsubdiffusive} for a similar result) shows that
  \begin{eqnarray} \label{eq:taur} \mathrm{T}_{\tau_{r}}^{(\kappa), \mathcal{L}} = \overline{B}_{r}( { \mathbf{T}}_{\kappa}).  \end{eqnarray}
Hence, by \eqref{eq:fixedtime} (and the remark following it), for any $t\in \mathcal{T}_{p}$ which is the hull of a ball of radius $r$ we have 
 \begin{eqnarray} \label{eq:stoppingtime} P\big( \overline{B}_{r}({ \mathbf{T}}_{ \kappa})=t\big) = P\big( \mathrm{T}_{\tau_{r}}^{(\kappa), \mathcal{L}} =t\big) = C_{p}^{(\kappa)} \cdot \kappa^{|t|}.  \end{eqnarray} Although the last display is sufficient to characterize the law of $ { \mathbf{T}}_{\kappa}$, it is not clear how it implies that $ { \mathbf{T}}_{\kappa}$ is $\kappa$-Markovian. To see this, consider yet another exploration process. Fix a triangulation $\Delta \in \mathcal{T}_{B}$ and let $t_{0} \subset t_{1} \subset ... \subset t_{n_{0}}= \Delta$ be an increasing sequence of triangulations with a boundary starting from the root edge so that $t_{i+1}$ is obtained from $t_{i}$ by the peeling of one (necessarily unique) edge $a_{i} \in \partial t_{i}$ for $i \leq n_{0}-1$. We consider the following modification of  the algorithm $ \mathcal{L}$:
$$ \mathcal{L}'(t) = \left\{\begin{array}{ll} a_{i} & \mbox{if }t=t_{i} \mbox{ for }i \in \{0,1, ... ,n_{0}-1\} \\ \mathcal{L}(t) & \mbox{otherwise.} \end{array} \right. $$ Here also, the peeling process with algorithm $ \mathcal{L}'$ will eventually swallow every vertex on the boundary of the growing triangulations and so the increasing union of $\mathrm{T}_{n}^{(\kappa), \mathcal{L}'}$ defines a random infinite triangulation of the plane denoted by $ \mathbf{T}_{\kappa}'$. We first show that $ \mathbf{T}'_{\kappa}$ has the same distribution as $ \mathbf{T}_{\kappa}$. First, remark that after step $n_{0}$, both peeling processes evolve according to the same rules. From this and a few simple geometric considerations we deduce that there exists some $r_{0}\geq 1$ (depending on $ \Delta$) such that for every $r \geq r_{0}$ we have $$ \mathrm{T}_{\tau'_{r}}^{(\kappa), \mathcal{L}'} = \overline{B}_{r}( \mathbf{T}'_{\kappa}),$$ where $\tau_{r}'$ is the first time at which no vertex of $\partial  \mathrm{T}_{n}^{(\kappa), \mathcal{L}'}$ is at distance less than or equal to $r-1$ from the origin.   Using  \eqref{eq:fixedtime} again we deduce that for any $t \in \mathcal{T}_{p}$ which is the hull of a ball of radius $r$ we have 
$P\big( \overline{B}_{r}( \mathbf{T}_{\kappa}')=t\big) = C_{p}^{(\kappa)} \kappa^{|t|}$. Comparing this with \eqref{eq:stoppingtime} we conclude that $ \overline{B}_{r}( \mathbf{T}'_{\kappa})$ and $ \overline{B}_{r}( \mathbf{T}_{\kappa})$ have the same law for every $r \geq r_{0}$. Since $ \mathbf{T}_{\kappa}$ and $ \mathbf{T}_{\kappa}'$  are both triangulations of the plane this entails that they have the same law. Coming back to the exploration process with algorithm $ \mathcal{L}'$, a moment's though shows that  $\Delta \subset \mathbf{T}'_{\kappa}$ if and only if we have $ \mathrm{T}_{i}^{(\kappa),  \mathcal{L}'} = t_{i}$ for every $i \in \{0,1,... , n_{0}\}$. In particular we have  
$$ P(  \Delta \subset \mathbf{T}'_{\kappa}) = P( \mathrm{T}_{n_{0}}^{(\kappa), \mathcal{L}'} = \Delta) \underset{\eqref{eq:fixedtime}}{=} C_{ | \partial \Delta| }^{(\kappa)} \cdot \kappa^{|\Delta|}.$$
Since $ \mathbf{T}_{\kappa} = \mathbf{T}'_{\kappa}$ in distribution, the last display still holds with $ \mathbf{T}'_{\kappa}$ replaced by $ \mathbf{T}_{\kappa}$. Because $\Delta$ was arbitrary this indeed shows that $ \mathbf{T}_{\kappa}$ fulfills \eqref{eq:SMP} and completes the construction.

\section{Geometric properties}
\subsection{Back to the peeling construction}
We first study in more details the peeling process of $ \mathbf{T}_{\kappa}$. This will be used in the proofs of Theorems \ref{thm:volume} and \ref{thm:hyperbolic}. Indeed, to study the volume growth of $ \mathbf{T}_{\kappa}$ we will explore it using the peeling by layers as in \cite{Ang03} and to establish Theorem \ref{thm:hyperbolic} one shall need to explore $ \mathbf{T}_{\kappa}$ along a simple random walk path as in \cite{BCsubdiffusive}. \medskip

 In the last section we \emph{constructed} $ \mathbf{T}_{\kappa}$ as the increasing union of triangulations given by an abstract peeling process.  In the rest of the paper, however, we will think of the peeling construction as ``embedded'' in $ \mathbf{T}_{\kappa}$ and exploring it. In other words, for every algorithm $ \mathcal{A}$, deterministic or using an extra source of randomness, we can couple a realization of $ \mathbf{T}_{\kappa}$ together with the sequence of growing triangulations $ ( \mathrm{T}_{n}^{(\kappa), \mathcal{A}})$ such that the latter is a growing subset of $ \mathbf{T}_{\kappa}$. Yet another way to express this is that we can explore $ \mathbf{T}_{\kappa}$ face by face (discovering the enclosed regions when needed) by peeling at each step an edge on the boundary of the current revealed part as long as this choice remains independent of the unexplored part. The proof of the above facts is merely a dynamical reformulation of Section \ref{sec:uniqueness} and is easily adapted from \cite[Section 1.2]{BCsubdiffusive} using the spatial Markov property of $ \mathbf{T}_{\kappa}$.  In particular the law of 
$$\big( P_{n}^{(\kappa)},V_{n}^{(\kappa)}\big)_{n \geq 0} := \big( | \partial	\mathrm{T}_{n}^{(\kappa), \mathcal{A}}|, | \mathrm{T}_{n}^{(\kappa), \mathcal{A}}|\big)_{n \geq 0}$$ does not depend on the algorithm $ \mathcal{A}$. More precisely, from Section \ref{sec:construction} we get that $P^{(\kappa)}$ is a Markov chain with transition probabilities given by 
 \begin{eqnarray} 
 \label{eq:loipeel1} P\big(\Delta P_{n}^{(\kappa)}=i \mid P_{n}^{(\kappa)}=p\big) &=& q_{i,p}^{(\kappa)} \qquad \mbox{ for } i \leq 1,  \end{eqnarray} where here and later $ \Delta X_{n}= X_{n+1}-X_{n}$. Conditionally on $P^{(\kappa)}$ the increments of $V^{(\kappa)}$ are independent and distributed as 
 \begin{eqnarray} \label{eq:loipeel2} \Delta V_{n}^{(\kappa)} & \overset{(d)}{=}& \mathcal{B}_{-\Delta P_{n}^{(\kappa)}},  \end{eqnarray} where $ \mathcal{B}_{i}$ is the law of the internal volume of a Boltzmann triangulation of the $i+1$-gon with $ \mathcal{B}_{-1}=1$ by convention.  Also, thanks to Lemma \ref{lem:fini} the increments of the chain $P_{n}^{(\kappa)}$ converges as the perimeter tends to $\infty$  towards i.i.d.\,\,steps of law $q_{i}^{(\kappa)} = \lim_{p \to \infty} q_{i,p}^{(\kappa)}$ which we recall is the step distribution of the random walk $\Xi^{(\kappa)}$ (started from $2$). Finally, conditionally on $\Xi^{(\kappa)}$ construct $\Omega^{(\kappa)}$ such that $\Delta \Omega^{(\kappa)}_{n}$ are independent and distributed as $\Delta \Omega^{(\kappa)}_{n} = \mathcal{B}_{- \Delta \Xi^{(\kappa)}_{n}}$ for every $n \geq 0$.

\begin{proposition}[Perimeter and volume growth during a peeling] \label{lem:peeling} Fix $ \kappa \in (0, \frac{2}{27})$ and recall the definitions of $\alpha$ and $\delta_{\kappa}$ in \eqref{eq:alphakappa} and  \eqref{eq:mean}. For every $ \varepsilon >0$ we have   \begin{eqnarray*} n^{1/2- \varepsilon} \left|\frac{P_{n}^{(\kappa)}}{n} - \delta_{\kappa}  \right| & \xrightarrow[n\to\infty]{a.s.}		& 0,\\
  n^{-1}{V_{n}^{(\kappa)}} & \xrightarrow[n\to\infty]{a.s.} &  \frac{\alpha(2 \alpha -1)}{ \delta_{\kappa}}.  \end{eqnarray*}
  \end{proposition}
  \proof
Recall the notation of Section \ref{sec:harmonic}. By \eqref{eq:htransform} and \eqref{eq:loipeel1} the Markov chain $(P_{n}^{(\kappa)})_{ n \geq 0}$ has the law of Doob's $h$-transform of the random walk $\Xi^{(\kappa)}$ (started from $2$)  of step distribution $ \boldsymbol{q}^{(\kappa)}$ by the function $ p \mapsto \tilde{C}^{(\kappa)}_{p}$ which is harmonic on $\{2,3,...\}$ and null for $p \leq 1$. By the results of \cite{BD94}, this process $P^{(\kappa)}$ has the same law as the walk $\Xi^{(\kappa)}$ conditioned on the event $\{ \Xi^{(\kappa)}_{i} \geq 2 : \forall i \geq 0\}$. Since $\Xi^{(\kappa)}$ has a positive drift $ \delta_{\kappa}$, the last event has a positive probability. Using  \eqref{eq:loipeel2} and the definition of the process $( \Xi^{(\kappa)}, \Omega^{(\kappa)})$ it follows that 

 \begin{eqnarray} \label{eq:staypositive} \big(P^{(\kappa)}_{n},V_{n}^{(\kappa)}\big)_{n \geq 0} \quad \overset{(d)}{=} \quad  \big(\Xi^{(\kappa)}_{n}, \Omega_{n}^{(\kappa)}\big)_{n \geq 0} \quad \mbox{ conditioned on } \underbrace{\{ \Xi^{(\kappa)}_{i} \geq 2, \forall i \geq 0\}}_{ \mbox{ positive proba.}}.  \end{eqnarray} 
In particular $ P^{(\kappa)}$ and $ \Xi^{(\kappa)}$ share the same almost sure properties. Since the step distribution of $ \Xi^{(\kappa)}$ has exponential tails, easy moderate deviations estimates (see e.g.\,\,\cite[Lemma 1.12]{LG05}) show that for every $ \varepsilon>0$ we have $\lim_{n \geq 0}n^{-1/2- \varepsilon}|\Xi^{(\kappa)}_{n} - n\cdot \delta_{\kappa}|=0$ almost surely. This implies the first statement of the proposition. 

For the second statement, we use the same argument. Using \eqref{eq:staypositive} it suffices to prove the similar result when $V_{n}^{(\kappa)}$ is replaced by $ \Omega_{n}^{(\kappa)}$. By the law of large numbers, this reduces to computing the mean of the increment of the random walk $\Omega^{(\kappa)}$.  From \cite[Proof of Proposition 3.4]{Ray13}\footnote{with the notation in \cite[Proposition 3.4]{Ray13}, we have $\theta = (1- \alpha)/2$ where $\alpha$ is given by \eqref{eq:alphakappa}.} we read that for $i \leq 0$
$$ E[ \mathcal{B}_{-i}] = \frac{i(2i-1)(1- \alpha)}{(3 \alpha-2) i +1}.$$
Plugging this into the definition of $\Omega^{(\kappa)}$ and using the explicit expression of the $q^{(\kappa)}_{\cdot}$ given by \eqref{eq:exact} it follows after a few manipulations using the generating function of Catalan numbers that $$E[\Delta \Omega_{n}] =  \alpha + \sum_{i \geq 1} i (2i-1)(1-\alpha)\frac{2}{4^i} \frac{(2i-2)!}{(i-1)!(i+1)!} \left( \frac{2}{ \alpha}-2\right)^i = \frac{\alpha(2\alpha-1)}{ \delta_{\kappa}}.$$ \endproof

\subsection{Volume growth} 
The goal of this section is to prove Theorem \ref{thm:volume}. We suppose that we discover $ \mathbf{T}_{\kappa}$ using the peeling algorithm with procedure $ \mathcal{L}$ which ``turns'' around the successive boundaries  $ \partial \overline{B}_{r}( \mathbf{T}_{\kappa})$ for $r \geq 0$ in a cyclic fashion, see Fig.\,\ref{fig:layers}.


\proof[Proof of Theorem \ref{thm:volume}] Recall that the stopping time $\tau_{r}$ is the first time in the exploration process when no vertex on the boundary is at distance less than or equal to $r-1$ from the origin of  $ \mathbf{T}_{\kappa}$ and that $ \mathrm{T}_{ \tau_{r}}^{(\kappa), \mathcal{L}} = \overline{B}_{r}( \mathbf{T}_{\kappa})$ by \eqref{eq:taur}. Recall also that $P_{n}^{(\kappa)}$ and $V_{n}^{(\kappa)}$ respectively are the perimeter and the size of the explored triangulation after $n$ steps of peeling. The proof is based on the following estimate:

\begin{lemma}[Time to complete a layer] \label{lem:tour} For any $ \varepsilon>0$ we have 
$$ \limsup_{r \to \infty} \big({\tau_{r+1}}\big)^{1/2- \varepsilon}\left| \frac{\tau_{r+1}-\tau_{r}}{P^{(\kappa)}_{\tau_{r}}} - \frac{2}{ \alpha - \delta_{\kappa}} \right| = 0.$$
\end{lemma}
Given the last lemma, the proof of Theorem \ref{thm:volume} is easy to complete. Indeed we have 
  \begin{eqnarray*}| \partial \overline{B}_{r+1}( \mathbf{T}_{\kappa})|- | \partial \overline{B}_{r}( \mathbf{T}_{\kappa})| &  \underset{ \eqref{eq:taur}}{=} & P_{\tau_{r+1}}^{(\kappa)}-P_{\tau_{r}}^{(\kappa)}\\  &\underset{\mathrm{Prop.\ } \ref{lem:peeling}}{=}&  \delta_{\kappa} |\tau_{r+1}- \tau_{r}| + o\big( (\tau_{r+1}^{(\kappa)})^{ \varepsilon+1/2}\big)\\
& \underset{ \mathrm{Lemma\ } \ref{lem:tour}}{=}  &  \frac{2 \delta_{\kappa}}{\alpha - \delta_{\kappa}}  P_{\tau_{r}} +o\big( (\tau_{r+1}^{(\kappa)})^{ \varepsilon+1/2}\big)\\
&\underset{ \eqref{eq:taur} \mathrm{\ and \ Prop.\ } \ref{lem:peeling} }=& \frac{2 \delta_{\kappa}}{\alpha - \delta_{\kappa}} | \partial \overline{B}_{r}( \mathbf{T}_{\kappa})|+ o\big(  | \partial \overline{B}_{r}( \mathbf{T}_{\kappa})|^{ \varepsilon+1/2}\big), \end{eqnarray*}  
equivalently 
$$  \left( \frac{\alpha - \delta_{\kappa}}{\alpha + \delta_{\kappa}} \right)\frac{| \partial \overline{B}_{r+1}( \mathbf{T}_{\kappa})|}{| \partial \overline{B}_{r}( \mathbf{T}_{\kappa})| } = 1 + o\big(  | \partial \overline{B}_{r+1}( \mathbf{T}_{\kappa})|^{ \varepsilon-1/2}\big),$$  where the $o$ is almost sure. Since $ P_{n}^{(\kappa)} \to \infty$ as $n \to \infty$, it  follows that $ | \partial \overline{B}_{r}( \mathbf{T}_{\kappa})| \to \infty$ as $r \to \infty$. Using the last display  we deduce that $\liminf_{r \to \infty} | \partial \overline{B}_{r}( \mathbf{T}_{\kappa})|^{1/r} \geq \frac{\alpha+\delta_{\kappa}}{\alpha - \delta_{\kappa}}$. Bootstrapping the argument by plugging this back into the last display, we deduce that the series $\sum_{r} \lambda_{r}-1$ is absolutely converging a.s.~where $$\lambda_{r} =  \left( \frac{\alpha - \delta_{\kappa}}{\alpha + \delta_{\kappa}} \right)\frac{| \partial \overline{B}_{r+1}( \mathbf{T}_{\kappa})|}{| \partial \overline{B}_{r}( \mathbf{T}_{\kappa})| } $$
A classic result then implies that $\prod_{r\geq 1} \lambda_{r}$ is converging in $ \mathbb{R}_{+}^*$ a.s.~otherwise said that we have the almost sure convergence
$$ \left(\frac{\alpha - \delta_{\kappa}}{\alpha+ \delta_{\kappa}}\right)^{r} \partial \overline{B}_{r}( \mathbf{T}_{\kappa}) \quad \xrightarrow[r\to\infty]{a.s.} \quad \Pi_{\kappa} \in (0, \infty).$$ 
This proves the first part of Theorem \ref{thm:volume}. The second part follows from Proposition \ref{lem:peeling} which shows that $V_{n}^{(\kappa)}/P_{n}^{(\kappa)} \to \frac{\alpha(2 \alpha-1)}{\delta_{ \kappa}^2}$ almost surely as $n \to \infty$.\endproof

\proof[Proof of Lemma \ref{lem:tour}] 
We adapt an argument from \cite{CLGpeeling}. Fix $r \geq 0$ and consider the situation at time $\tau_{r}$. In the future of the peeling process we will go cyclically around $ \partial \overline{B}_{r}( \mathbf{T}_{\kappa}) = \partial \mathrm{T}_{\tau_{r}}^{(\kappa), \mathcal{L}}$  from left to right swallowing the vertices of $\partial \overline{B}_{r}( \mathbf{T}_{\kappa})$ (see Fig.\,\ref{fig:layers}) until none is left on the active boundary which happens at time $\tau_{r+1}$. For $\tau_{r} \leq i \leq \tau_{r+1}$ denote by $A_{i}$ the number of vertices of $\partial \overline{B}_{r}( \mathbf{T}_{\kappa})$ which are still part of the boundary of $ \mathrm{T}_{i}^{(\kappa), \mathcal{L}}$, so that $A_{\tau_{r}} = P_{\tau_{r}}$ and $A_{\tau_{r+1}}=0$. Clearly we have 
 \begin{eqnarray} \label{eq:taur1} \tau_{r+1} - \tau_{r} =\inf\left\{ i \geq 0 : \sum_{j=0}^{i-1} \Delta A_{j+ \tau_{r}} = - P_{\tau_{r}}\right\}. \end{eqnarray} Also introduce the events $L_{i},R_{i},C_{i}$ respectively realized when the peeling process at time $i$ discovers a triangle bent to the left, bent to the right, or pointing inside the undiscovered part. We claim that a good approximation of the behavior of the process $A$ is given by
 \begin{eqnarray} \Delta A_{i} &\approx& \Delta P^{(\kappa)}_{i} \mathbf{1}_{D_{i}}, \qquad \tau_{r} \leq i < \tau_{r+1} \label{eq:firstapprox} \end{eqnarray}
Let us first imagine that the last display holds exactly and let us show why this implies the lemma, we then sketch how to cope with the approximation. We claim that almost surely we have 
$$ \sum_{i=\tau_{r}}^{\tau_{r}+n} \Delta P_{i}^{(\kappa)} \mathbf{1}_{D_{i}} = -\frac{\alpha- \delta_{\kappa}}{2}n + o( n^{1/2+ \varepsilon}).$$ To prove the claim, we use the same argument as the one that led to \eqref{eq:staypositive} and argue that it is sufficient to prove the last display when $ (\Delta P_{i}^{(\kappa)} \mathbf{1}_{D_{i}})$ is replaced by $ (\Delta \Xi^{(\kappa)}_{i}  \mathbf{1}_{ \Delta \Xi^{(\kappa)}_{i} < 0} \epsilon_{i})$ where $ \epsilon_{i}$  are i.i.d.\,\,Bernoulli variables $P( \epsilon_{i}=1)=P( \epsilon_{i}=0)=1/2$ also independent of $\Xi^{(\kappa)}$. An easy calculation then shows that $E[\Delta \Xi^{(\kappa)}_{i}  \mathbf{1}_{ \Delta \Xi^{(\kappa)}_{i} < 0} \epsilon_{i}] = \frac{\alpha-\delta_{\kappa}}{2}$ and moderate deviations arguments (see e.g.\,\,\cite[Lemma 1.12]{LG05}) imply  the claim. From this and \eqref{eq:taur1}, we deduce that $\tau_{r+1} = \tau_{r} + \frac{2}{\alpha-\delta_{\kappa}} P_{\tau_{r}} + o({\tau_{r+1}}^{1/2+ \varepsilon})$ as wanted.

\begin{figure}[!h]
 \begin{center}
 \includegraphics[width=12cm]{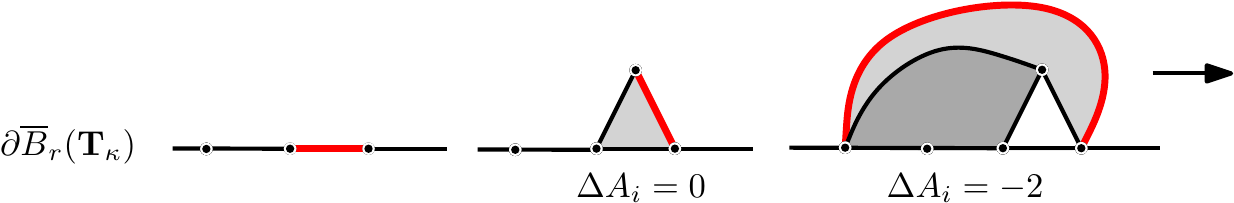}
 \caption{ \label{fig:left} In the beginning of the peeling of the $r$th layer, a few peeling steps towards the left may contribute to swallowing the vertices of $ \partial \overline{B}_{r}( \mathbf{T}_{\kappa})$.}
 \end{center}
 \end{figure}

Let us now see why the approximation \eqref{eq:firstapprox} is quite good. See Fig.\,\,\ref{fig:layers}. Indeed, they are only two cases when \eqref{eq:firstapprox} may fail. First of all, during the last peeling step $i= \tau_{r+1}-1$ we could have $\Delta P^{(\kappa)}_{i}  \mathbf{1}_{D_{i}} < \Delta A_{i}$ since this last jump  towards the right could swallow more than just the remaining edges of $ \partial \overline{B}_{r}( \mathbf{T}_{\kappa})$. This is not a big problem since it concerns only one step (and $\Delta P^{(\kappa)}$ has exponential tails). However, a bit more annoying is the fact that peeling steps corresponding to jumps towards the left could actually contribute to reducing $A$ as well, see Fig.\,\ref{fig:left}.

However, an easy calculation shows that at step $\tau_{r} + i$ there are roughly $ \frac{\alpha + \delta_{\kappa}}{2} \cdot i$ edges on $ \partial \mathrm{T}_{i+ \tau_{r}}^{(\kappa), \mathcal{L}}$ separating   $  \partial \overline{B}_{r}( \mathbf{T}_{\kappa})$ from the left of the current edge to peel. Since $\Delta P_{i}^{(\kappa)}$ has exponential tails, we deduce that the last phenomenon can only appear in the first few $\ln(\tau_{r})$ steps after $\tau_{r}$ and cannot perturb the approximation \eqref{eq:firstapprox} too much. We leave the details to the careful reader.  \endproof

\subsection{Anchored expansion}
Like in many stochastic examples which are hyperbolic in flavor (e.g.\,\,supercritical Galton--Watson trees), the randomness of $ \mathbf{T}_{\kappa}$ allows for any particular pattern to happen somewhere in the lattice and thus destroys any hope of having a positive Cheeger expansion constant. The latter has to be replaced by a more refine notion: the anchored expansion constant. \medskip 

If $G$ is a connected graph with an origin vertex $\rho$ and if $S$ is a subset of vertices of $G$, we denote by $ | \partial_{E}S|$ the number of edges having an endpoint in $S$ and the other outside $S$. Also, write $|S|_{E}$ for the sum of the degrees of the vertices of $S$. The edge anchored expansion constant of $G$ is defined by
$$ i^*_{E}(G) = \liminf_{n \to \infty} \left \{ \frac{| \partial_{E} S|}{|S|_{E}} : S \subset \mathrm{Vertices}(G), S \mbox{ finite and connected}, \rho \in S, |S|_{E} \geq n\right\}.$$
It is easy to see that the above definition does not depend on the origin point $\rho$. See \cite{Vir00} for background on anchored expansion. As in \cite[Theorem 2.2]{Ray13} the spatial Markov property enables us to deduce almost effortless that our lattices have a positive anchored expansion constant in the hyperbolic regime:
\begin{proposition}[Edge anchored expansion] \label{prop:anchored}For $ \kappa < \frac{2}{27}$ we have $i_{E}^*( \mathbf{T}_{\kappa}) >0$ almost surely.  \end{proposition}
By ergodicity (see the proof of Proposition \ref{prop:ergodic} below) the variable $i_{E}^*( \mathbf{T}_{\kappa})$ is actually almost surely constant. We do not have a good guess for its correct value. Before doing the proof of Proposition \ref{prop:anchored} we state a lemma. For $p, n \geq 0$ we denote by $ \mathcal{T}_{n,p}$ the set of all finite triangulations of the $p$-gon with $n$ vertices.

\begin{lemma} \label{lem:estima}There exists $m_{0} \geq 1$ and $c_{1}>0$ such that for every $m \geq m_{0}$ and every $ n  \geq  m p$ we have 
$$ \# \mathcal{T}_{n+p,p} \leq c_{1} \sqrt{\frac{p}{n^3}}  \ 9^p \left(\frac{27}{2} \right)^n.$$
\end{lemma}
\proof For $p \geq 2$ and $n \geq 0$, if $ \mathcal{T}_{n+p,p}^\rightarrow$ is the set of all $2$-connected  triangulations of the $p$-gon of size $n+p$ such that the hole is on the right-hand side of the root edge then from \cite{GJ83} we read 
$$ \# \mathcal{T}_{n+p,p}^\rightarrow = \frac{(2p-3)!}{(p-2)!(p-2)!} 2^{n+1} \frac{(2p+3n-4)!}{n!(2p+2n-2)!}.$$
An application of Euler's formula shows that such a triangulation has exactly $3n+2p-3$ edges. Hence we deduce that 
$ \# \mathcal{T}_{n+p,p} \leq 2(3n+2p-3) \# \mathcal{T}_{n+p,p}^\rightarrow$. Suppose now that  $n \geq m p$ for $m \geq 1$. Using the last display and Stirling's formula we have  for a constant $c>0$ that may vary from line to line
  \begin{eqnarray*} \# \mathcal{T}_{n+p,p}  \leq2(3n+2p-3) \# \mathcal{T}_{n+p,p}^\rightarrow  &\leq& c n \frac{(2p-3)!}{(p-2)!(p-2)!} 2^{n+1} \frac{(2p+3n-4)!}{n!(2p+2n-2)!} \\ & \leq &c n \sqrt{p} 4^p 2^{n} n^{-1/2} \frac{(2p+3n-4)^{2p+3n-4}}{n^n (2p+2n-2)^{2p+2n-2}} \\ & \leq & c  \sqrt{ \frac{p}{n^3}} 9^p \left( \frac{27}{2}\right) ^n \frac{(1 + \frac{2p-4}{3n})^{2p+3n-4}}{( 1+ \frac{p-1}{n})^{2p+2n-2}}\\
&  \leq & 
  c  \sqrt{ \frac{p}{n^3}} 9^p \left( \frac{27}{2}\right) ^n \frac{(1 + \frac{2p}{3n})^{2p+3n}}{( 1+ \frac{p}{n})^{2p+2n}}. \end{eqnarray*}
  If $n \geq m p$ with $m$ sufficiently large, the last fraction is the preceding display is smaller than one. This completes the proof of the lemma.
\endproof

\proof[Proof of Proposition \ref{prop:anchored}] First, in the definition of $i_{E}^*( \mathbf{T}_{\kappa})$ we can restrict ourself to those sets $S$ such that $ \mathbf{T}_{\kappa} \backslash S$ has only one (infinite) component because filling-in the finite holes decreases the boundary size and increases the volume. We then consider the triangulation with one hole $ \overline{S}$ obtained by adding all the faces adjacent to a vertex of $S$ as well as the finite regions enclosed. One may check that $ |\partial \overline{S}| \leq | \partial_{E}S|$. By Euler's relation we also get that $3| \overline{S}| = | \partial \overline{S}| + 3 + \# \mathrm{Edges}( \overline{S})$, hence $ |\overline{S}| \geq \# \mathrm{Edges}( \overline{S})/3$. Since we also have $ \# \mathrm{Edges}( \overline{S}) \geq \frac{1}{2}|S|_{E}$ we get that $| \overline{S}| \geq \frac{1}{6} |S|_{E}$ and so
$$ \frac{| \partial \overline{S}|}{| \overline{S}|} \leq  6 \frac{| \partial_{E}S|}{|S|_{E}}.$$
To prove the proposition, it is thus sufficient to show that the ratio $|\partial A|/|A|$ is bounded away from $0$ for all triangulations $A \in \mathcal{T}_{B}$ such that $A \subset \mathbf{T}_{\kappa}$. For this we crudely use a first moment method. Fix $m \geq 1$. We have
  \begin{eqnarray*}  P\Big( \exists A \subset \mathbf{T}_{\kappa} : |A| > (m+1) | \partial A| \Big ) &\leq&  E\Big[ \# A \subset \mathbf{T}_{\kappa} : |A| > (m+1) | \partial A| \Big  ] \\ &=& \sum_{p \geq 1} \sum_{ n > m p} \sum_{ A \in \mathcal{T}_{n+p,p}} P( A \subset \mathbf{T}_{ \kappa}) \\
  &\underset{ \eqref{eq:SMP}}{=}& \sum_{p \geq 1} \sum_{ n > m p} C_{p}^{(\kappa)} \kappa^{n+p} \# \mathcal{T}_{n+p,p}.
  \end{eqnarray*}
At this point we use Lemma \ref{lem:estima} and get for $n \geq mp$ with $m \geq m_{0}$
$$ P\Big( \exists A \subset \mathbf{T}_{\kappa} : |A| > (m+1) | \partial A| \Big ) \leq c_{1} \sum_{p \geq 1} C_{p}^{(\kappa)} (9\kappa) ^p \sqrt{p}\sum_{ n > m p} (\frac{27}{2} \cdot \kappa)^n n^{-3/2}.$$
Since $\kappa < \frac{2}{27}$ the last sum is easily seen to be smaller than $ c_{2} ( \frac{27}{2} \cdot \kappa)^{mp}$ for some constant $c_{2}>0$ depending on $\kappa$. Also, from Lemma \ref{lem:fini} we have $C_{p}^{(\kappa)} \leq c_{3} \cdot \beta^{-p}$ for some constant $c_{3}>0$ still depending on $\kappa$. Hence we have 
 \begin{eqnarray*} P\Big( \exists A \subset \mathbf{T}_{\kappa} : |A| > (m+1) | \partial A| \Big ) &\leq& c_{1}c_{2}c_{3} \sum_{p \geq 1}  \sqrt{p} \, \left( \frac{9 \kappa}{\beta} \cdot \left(\kappa \frac{27}{2}\right)^m \right)^{p}.   \end{eqnarray*}
Since $\kappa <  \frac{2}{27}$, by choosing $m$ large enough, we can make $\frac{9 \kappa}{\beta} \cdot \left(\kappa \frac{27}{2}\right)^m$ as small as we wish and thus the last probability tends to $0$ as $m \to \infty$.  This indeed implies that $P( i_{E}^* ( \mathbf{T}_{\kappa}) =0)=0$ and completes the proof of the proposition. \endproof

 \section{Simple random walk}
 
 \label{sec:speed}

 In this section, we study the simple random walk on $ \mathbf{T}_{\kappa}$. The special case $\kappa = \frac{2}{27}$ has already received a lot of attention and it is known that the UIPT is recurrent \cite{GGN12} and subdiffusive in the quadrangular case \cite{BCsubdiffusive}. In this section we prove Theorem \ref{thm:hyperbolic}. \medskip 
 
First of all, Proposition \ref{prop:anchored} combined with the result of \cite{Tho92} shows  that   \begin{eqnarray} \label{eq:transient} \mathbf{T}_{\kappa} \mbox{ is almost surely transient for } \kappa < \frac{2}{27}. \end{eqnarray} In the bounded degree case, a positive anchored expansion constant is even sufficient to imply positive speed for the simple random walk as shown by Virag \cite{Vir00}. Unfortunately, the lack of a uniform bound on the degrees in $ \mathbf{T}_{\kappa}$ prevents us from using this nice result and we shall go through a rather winding but bucolic bypass. The  strategy to prove Theorem \ref{thm:hyperbolic} is  the following:
 \begin{eqnarray} \label{sketch} 
 \mbox{study of the peeling along a SRW}   & \underset{ \mbox{Section \ref{sec:noninter}}}{\Longrightarrow}  & \mbox{non (intersection property)} \\
     & \underset{\cite{BCGLiouville}}{\Longrightarrow} & \mbox{non Liouville} \nonumber \\
      &\underset{ \mbox{Section \ref{sec:entropy}}}{ \Longrightarrow} &\mbox{positive speed}.  \nonumber\end{eqnarray}

 Let us introduce a piece of notation. Conditionally on $ \mathbf{T}_{\kappa}$ consider a simple random walk (at each step, independently of the past, walk through one adjacent oriented edge uniformly at random)  started from the \emph{target} of the root edge and denote by $(\vec{E}_{i})_{i \geq 0}$ the sequence of oriented edges traversed by the walk where by convention $\vec{E}_{0}$ is the root edge. We also denote by $X_{0}, ... , X_{n}$ the successive vertices visited by the walk, i.e.\,\,$X_{i}$ is the origin of the oriented edge $\vec{E}_{i}$. The underlying probability and expectation relative to the lattice $ \mathbf{T}_{\kappa}$ are denoted by $P$ and $E$ whereas  the (quenched) probability and expectation relative to the walk on $ \mathbf{T}_{\kappa}$ are denoted by $ \mathbb{P}$ and $ \mathbb{E}$.

\subsection{Reversibility and ergodicity}
The notation $\overleftarrow{e} $ stands for the reversed oriented edge $\overrightarrow{e}$.

\begin{proposition}[Reversibility]  \label{prop:reversible} For any $\kappa \in (0, \frac{2}{27}]$ and every $i \geq 0$ we have the equalities  in distribution 
$$\begin{array}{crcl} (i) \qquad & ( \mathbf{T}_{\kappa} ; \overrightarrow{E}_{0}) &=& ( \mathbf{T}_{\kappa} ; \overleftarrow{E}_{0})\\ (ii)\qquad  & (  \mathbf{T}_{\kappa} ;  \overrightarrow{E}_{0}, ... , \overrightarrow{E}_{i}) &=&  (  \mathbf{T}_{\kappa} ;  \overleftarrow{E}_{i},... , \overleftarrow{E}_{0}).\end{array}$$
\end{proposition}
Combining the statements of the last proposition we deduce that $( \mathbf{T}_{\kappa} ; \overrightarrow{E}_{0}) = ( \mathbf{T}_{\kappa} ; \overleftarrow{E}_{i}) = ( \mathbf{T}_{\kappa} ; \overrightarrow{E}_{i})$ in distribution. This proves that the law of the lattice is unchanged under re-rooting along a simple random walk path. We say in short that $ \mathbf{T}_{ \kappa}$ is a stationary (in our case also reversible) random graph, see  \cite[Definition 1.3]{BCstationary}. We refer to \cite[Section 2.1]{BCstationary} for more details about the connections between the concepts of stationary (and reversible) random graphs, ergodic theory, unimodularity, mass-transport principle and measured equivalence relations. Note that in the critical case $\kappa = \frac{2}{27}$, the stationarity of the UIPT is an easy  consequence of the fact that it is a local limit of uniformly rooted finite graphs (see \cite[Theorem 3.2]{AS03}). Although we conjecture that $ \mathbf{T}_{\kappa}$ can similarly be obtained as the local limit of uniformly rooted random triangulations in high genus (Conjecture 1) we provide a direct proof of Proposition \ref{prop:reversible}.

\proof Point $(i)$ is easy: the lattice obtained from $ \mathbf{T}_{\kappa}$ is still $\kappa$-Markovian and thus has the same distribution by Theorem \ref{thm:un}. Let us now turn to $(ii)$. Let $i,r >0$. Fix a triangulation with a boundary $t \subset \mathcal{T}_{B}$ and a path $w=( \vec{e}_{0}, \vec{e}_{1}, ... , \vec{e}_{i})$ such that $w$ could be the result of a $i$-step random walk inside $t$ (with the convention that $ \vec{e}_{0}$ is the root edge of $t$). We denote by $x_{0}, x_{1},... , x_{i+1}$ the vertices visited by the path. Furthermore, we assume that $t$ is the hull of the ball of radius $r$ around $w$ in the sense that it is made of all the faces containing a vertex at graph distance (inside $t$) smaller than or equal to $r-1$ from the set $\{x_{0}, x_{1}, ... , x_{i+1}\}$ as well as the finite regions enclosed. We write $t= \overline{B}_{r}(\{x_{0}, ... , x_{i+1}\})$.  We now ask what is the probability that, inside $ \mathbf{T}_{\kappa}$, that the first $i$ steps of the walk correspond to $w$ and that the hull of radius $r$ around these is $t$:
 \begin{eqnarray*}P( \overrightarrow{E}_{k}= \vec{e}_{k}, \forall	k \leq i \mbox{ and } \overline{B}_{r}(\{x_{0}, ... , x_{i+1}\})=t)
 &=& P( t \subset \mathbf{T}_{\kappa}) \cdot P( \overrightarrow{E}_{k}= \vec{e}_{k}, \forall	k \leq i \mid t \subset \mathbf{T}_{\kappa}) \\
 & \underset{ \eqref{eq:SMP} }{=}	 & C_{| \partial t|}^{(\kappa)} \kappa^{|t|} \prod_{k=1}^{i}  \mathrm{deg}(x_{k})^{-1}.  \end{eqnarray*} 
We now remark that the last probability is exactly the same if we replace $(t,w)$ with the same triangulation $t$ and the reversed path $\overleftarrow{w} = (\overleftarrow{e_{i}}, ... , \overleftarrow{e_{0}})$. Since $r$ is arbitrary this proves that $(  \mathbf{T}_{\kappa} ;  \overrightarrow{E}_{0}, ... , \overrightarrow{E}_{i})$ and $ (  \mathbf{T}_{\kappa} ;  \overleftarrow{E}_{i},... , \overleftarrow{E}_{0})$ indeed have the same law.  \endproof

\begin{proposition}[Ergodicity] \label{prop:ergodic} The shift operation $\theta : ( \mathbf{T}_{\kappa} ; ( \overrightarrow{E_{i}})_{i \geq  \mathbf{0}})\mapsto ( \mathbf{T}_{\kappa} ; ( \overrightarrow{E_{i}})_{i \geq   \mathbf{1}})$  is ergodic for $ P \otimes \mathbb{P}$. 
\end{proposition}
\proof We consider the set $ (\mathcal{G}_{*}, \mathrm{d_{loc}}, \mathcal{F})$ of all locally finite connected rooted graphs endowed with the local distance and the associated Borel $\sigma$-field. An easy adaptation of \cite[Theorem 7.2]{Ang03} shows that the class $ \mathcal{C}\subset \mathcal{F}$ of events which are invariant (up to $P$-measure zero) to finite changes in the triangulation is trivial for $ P$, i.e.\,\, $$A \in \mathcal{C} \Rightarrow P(A)  \in \{0,1\}.$$
We now adapt  \cite[Theorem 4.6]{AL07} and prove that the last display implies ergodicity of the shift along a simple random walk. Formally, consider $ (\mathcal{P}_{*}, \mathrm{d_{loc}}, \mathcal{K})$ the set of all locally finite connected rooted graphs together with an infinite path on them endowed with (an extension) of the local distance and the associated Borel $\sigma$-field. Let $ B \in \mathcal{K}$ be an event invariant by the shift along the path. As in the proof of \cite[Theorem 5.1]{LS99} we have $ \mathbb{P}_{ \mathbf{T}_{\kappa}}(B) \in \{0,1\}$ where $ \mathbb{P}_{G}$ is the probability measure induced on $ \mathcal{P}_{*}$ by the simple random walk on the (fixed) graph $G$. Consider then the event $ A =\{ \mathbb{P}_{ \mathbf{T}_{\kappa}}( \mathcal{B})=1\} \in \mathcal{F}$. Since $ \mathbf{T}_{\kappa}$ is almost surely transient \eqref{eq:transient}, a moment's thought shows that $ A$ is invariant (up to event of $P$-measure $0$) by finite changes in the triangulation. It follows by the last display that $P({A}) \in \{0,1\}$ whence $ P \otimes \mathbb{P}( {B}) \in \{0,1\}$ as desired. \endproof

Let us give an application of the last result and show existence of the speed (but not the positivity of the latter). Combining the stationarity of $ \mathbf{T}_{\kappa}$ given after Proposition \ref{prop:reversible} together with Proposition \ref{prop:ergodic}, an application of Kingman's ergodic subadditive theorem (see e.g.  \cite[Theorem 2.2]{BCstationary} or \cite[Proposition 4.8]{AL07}) proves the following convergence 
 \begin{eqnarray} \label{eq:abstractspeed}  n^{-1} \mathrm{d_{gr}}(X_{0},X_{n}) \quad \xrightarrow[n\to\infty]{P \otimes \mathbb{P} \ a.s.} \quad s_{\kappa} \in [0,1]. \end{eqnarray}

\subsection{Non-intersection by peeling}
\label{sec:noninter}

Following the proof-sketch \eqref{sketch} we start by studying the intersection properties of $ \mathbf{T}_{\kappa}$. Recall that a graph $G$ is said to have the intersection property if \emph{almost surely} the range of two independent simple random walks intersect infinitely often. It is easy to see that this property does not depend on the starting points of the walks.  In this section we show:
\begin{proposition}[non-intersection] \label{lem:intersection} When $ \kappa \in (0, \frac{2}{27})$  almost surely $ \mathbf{T}_{\kappa}$ does not possess the intersection property.
\end{proposition}

The key tool to prove Proposition  \ref{lem:intersection} is the peeling process along a simple random path, specifically we explore $ \mathbf{T}_{\kappa}$ using a peeling algorithm that discovers the triangulation when  necessary for the walk to make one more step. This was first used in \cite{BCsubdiffusive} to establish the subdiffusivity of simple random walk on random quadrangulations. We start with the formal definition of this algorithm denoted $ \mathcal{W}$ (for ``walk'') and then interpret it in terms of pioneer points. Recall that by convention, the first step of the walk is $\vec{E}_{0} = (X_{0},X_{1})$ and so we shall start the process at the target of the root edge.\medskip 

We define a sequence $ \vec{e} = \mathrm{T}_{0}^{(\kappa), \mathcal{W}} \subset  \mathrm{T}_{1}^{(\kappa), \mathcal{W}} \subset \cdots \subset  \mathrm{T}_{n}^{ (\kappa), \mathcal{W}} \subset\cdots \subset  \mathbf{T}_{ \kappa}$ of triangulations with boundaries and two random non decreasing functions $f,g : \mathbb{N} \to \mathbb{N}$ such that $f(0)=0$, $g(0)=1$ and  \begin{eqnarray} \label{eq:dedans}X_{g(k)} \in  \mathrm{T}_{f(k)}^{ (\kappa),\mathcal{W}}, \quad \mbox{ for every } k\geq 0, \end{eqnarray}  whose evolution is described by induction as follows.  We have two cases: \begin{itemize}
\item If the current position $X_{g(k)}$ of the simple random walk belongs to $\partial  \mathrm{T}_{f(k)}^{ (\kappa),\mathcal{W}}$, then choose an edge $a$ on $\partial  \mathrm{T}_{f(k)}^{ (\kappa),\mathcal{W}}$ containing $X_{g(k)}$ and set $f(k+1):=f(k)+1$ and $g(k+1):=g(k).$ The triangulation $ \mathrm{T}_{f(k+1)}^{ (\kappa),\mathcal{W}}$ is the map obtained from $ \mathrm{T}_{f(k)}^{ (\kappa),\mathcal{W}}$ after peeling the edge $a$. \item If the current position $X_{g(k)}$ of the simple random walk belongs to $ \mathrm{T}_{f(k)}^{ (\kappa),\mathcal{W}}\backslash \partial  \mathrm{T}_{f(k)}^{ (\kappa),\mathcal{W}}$ then we set $f(k+1):=f(k)$ and $g(k+1):=g(k)+1$. In words, we let the walker move for one more step and do not touch the explored triangulation.
\end{itemize}

Note that we have $f(n)+g(n)=n+1$ and $f,g \to \infty$. Although this algorithm has an extra randomness due to the simple random walk, the edges chosen to be revealed in the peeling process are independent of the unknown part, and thus the process $(|\partial  \mathrm{T}_{n}^{ (\kappa),\mathcal{W}}|,| \mathrm{T}_{n}^{ (\kappa),\mathcal{W}}|)_{n\geq 0}$ has the same law as $(P_{n}^{(\kappa)},V_{n}^{(\kappa)})_{n \geq 0}$ of Proposition \ref{lem:peeling}. In the following, it will be important to have a geometric interpretation of this algorithm.

\paragraph{Interpretation.} For any $k \geq 0$ consider the submap $ \mathsf{Hull}(X_{1}, ... , X_{k}) \subset \mathbf{T}_{\kappa}$ formed by the faces that are adjacent to $\{X_{1}, X_{2}, \ldots , X_{k}\}$ as well as the finite holes they enclose. By convention $ \mathsf{Hull}( \varnothing)$ is the root edge. Then an easy geometric lemma (see \cite[Proposition 7]{BCsubdiffusive}) shows that the peeling times exactly correspond to the times when
\begin{eqnarray*} X_{g(k)} &\in& \partial  \mathsf{Hull}( X_{1}, ... , X_{g(k)-1}).  \end{eqnarray*}
These points are called pioneer points. In other words, as soon as the walk reaches a pioneer point, the peeling process starts to discover the neighborhood of the current position  (this typical takes a few steps of peeling) enabling the simple random walk to displace again.

\begin{figure}[!h]
 \begin{center}
 \includegraphics[width=7cm]{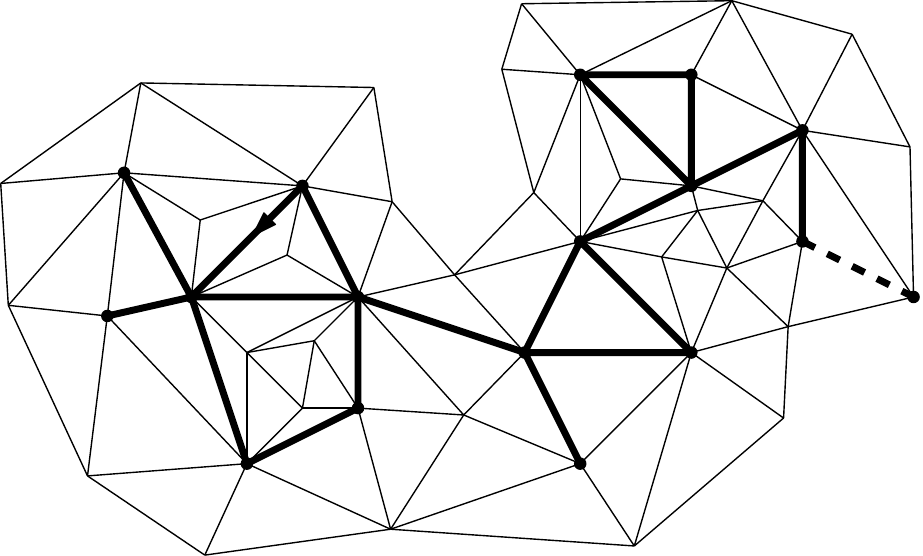}
 \caption{The trace of the simple random walk about to reach a pioneer point.}
 \end{center}
 \end{figure}

\begin{lemma} We have \label{lem:X0}
$ P \otimes \mathbb{P}\big( X_{0} \in \partial  \mathsf{Hull}(X_{1}, ... , X_{n}) : \forall n \geq 1 \big)>0.$ \label{lem:c>0}
\end{lemma}

\begin{remark} Note that in the peeling process along the simple random walk, the first pioneer point is $X_{1}$ and it is indeed  possible that $X_{0}$ stays on the boundary of the discovered triangulation for ever. The lemma says that this happens with positive probability.
\end{remark}

\proof Note that the events $ \{ X_{0} \in \partial \mathsf{Hull}(X_{1}, ... , X_{n})\}$ are clearly decreasing in $n$  so their $ P \otimes \mathbb{P}$-probabilities tend to some constant $ c \in [0,1]$. We have to show that $c>0$.  By the stationarity and reversibility of the walk on the lattice (Proposition \ref{prop:reversible}) we have 
 \begin{eqnarray*} P \otimes \mathbb{P}\big( X_{0} \in  \partial  \mathsf{Hull}(X_{1}, ... , X_{n})\big) &\underset{  \mathrm{reversibility}}{=}& P \otimes \mathbb{P}\big( X_{n} \in  \partial  \mathsf{Hull}(X_{n-1}, ... , X_{0})\big) \\ 
&\underset{ \mathrm{stationarity}}{=}& P \otimes \mathbb{P}\big( X_{n+1} \in  \partial  \mathsf{Hull}(X_{n}, ... , X_{1})\big) \\
&\underset{ \mathrm{definition}}{=}& P \otimes \mathbb{P}( X_{n+1} \mbox{ is pioneer})   \end{eqnarray*}
and so
 \begin{eqnarray} \label{eq:c}  \lim_{n \to \infty}^\downarrow P \otimes \mathbb{P}( X_{n} \mbox{ is pioneer}) = c.  \end{eqnarray}
We now combine the transience of the walk with the peeling estimates of Proposition \ref{lem:peeling}. Specifically, using the transience \eqref{eq:transient}, the stationarity (Proposition \ref{prop:reversible}) and the ergodicity of the simple random walk (Proposition \ref{prop:ergodic}), we deduce from \cite[Theorem 2.2]{BCstationary} that the range of the simple random walk grows linearly i.e.\,\,there exists $ \eta>0$ such that we have the almost sure convergence under $ P \otimes \mathbb{P}$
 \begin{eqnarray} n^{-1}\#\{ X_{0}, ... , X_{n}\} & \xrightarrow[n\to\infty]{a.s.}& \eta. \label{eq:rangelinear}  \end{eqnarray}
If we let run the peeling algorithm for $n$ steps (either peeling or walk step) then we have 
$$  \eta\,  g(n) \quad \underset{ \eqref{eq:rangelinear}}{\sim} \quad \# \{X_{0}, ... , X_{g(n)}\} \quad \underset{ \eqref{eq:dedans}}{\leq} \quad | \mathrm{T}_{f(n)}^{ (\kappa),\mathcal{W}}| \quad \underset{ \mathrm{Prop.}\, \ref{lem:peeling}}{\sim} \quad \frac{\alpha(2 \alpha-1)}{ \delta_{\kappa}} f(n).$$
Since $f(n)+g(n) =n+1$ we have   \begin{eqnarray} \liminf_{n \to \infty}\label{eq:f} \frac{f(n)}{n} \geq  \left(1+ \frac{\alpha(2\alpha-1)}{\eta   \delta_{\kappa}}\right)^{-1}. \end{eqnarray} Notice that the discovery of a pioneer point automatically triggers at least one peeling step. On the other hand, an estimate similar to \cite[Proposition 3.6]{AR13} or \cite[Lemma 4.4]{Ray13} shows that for any $k \geq 0$, when discovering the $k$th pioneer point, the expected number of peeling steps needed to perform a new random walk step is stochastically dominated by a geometric variable with a fixed parameter. It follows that if we put $p(n) =\#\{i \leq g(n) : X_{i} \mbox{ is pioneer}\}$ then  for some constant $\Lambda \geq1$ we almost surely have  \begin{eqnarray} \label{eq:pn}  \limsup_{n \to \infty } \frac{f(n)}{p(n)} \leq \Lambda .  \end{eqnarray}
We deduce that a.s.\,\,the asymptotic proportion of random walk steps which are pioneer satisfies 
$$ \liminf_{n \to \infty} \frac{p(n)}{g(n)} \geq \liminf_{n \to \infty} \frac{p(n)}{n} \underset{ \eqref{eq:pn}}{\geq} \Lambda^{-1} \liminf_{n \to \infty} \frac{f(n)}{n} \underset{ \eqref{eq:f}}{\geq}   \Lambda^{-1}\left(1+ \frac{\alpha(2\alpha-1)}{r  \delta_{\kappa}}\right)^{-1}.$$ By Cesàro theorem and \eqref{eq:c} we have 
 \begin{eqnarray*} c &  \underset{ \eqref{eq:c}}{=} & \liminf_{n \to \infty} \frac{1}{n} E \otimes \mathbb{E}\left[ \sum_{k=1}^{n}  \mathbf{1}\{X_{k} \mbox{ is pioneer}\} \right] \\ & \underset{ \mathrm{Fatou}}{\geq} & E \otimes \mathbb{E}\left[ \liminf_{n \to \infty} \frac{p(n)}{g(n)} \right] \geq  \Lambda^{-1}\left(1+ \frac{\alpha(2\alpha-1)}{r  \delta_{\kappa}}\right)^{-1},   \end{eqnarray*} by the last display. This completes the proof of the lemma. 
\endproof 

\begin{figure}[!h]
  \begin{center}
  \includegraphics[width=16cm]{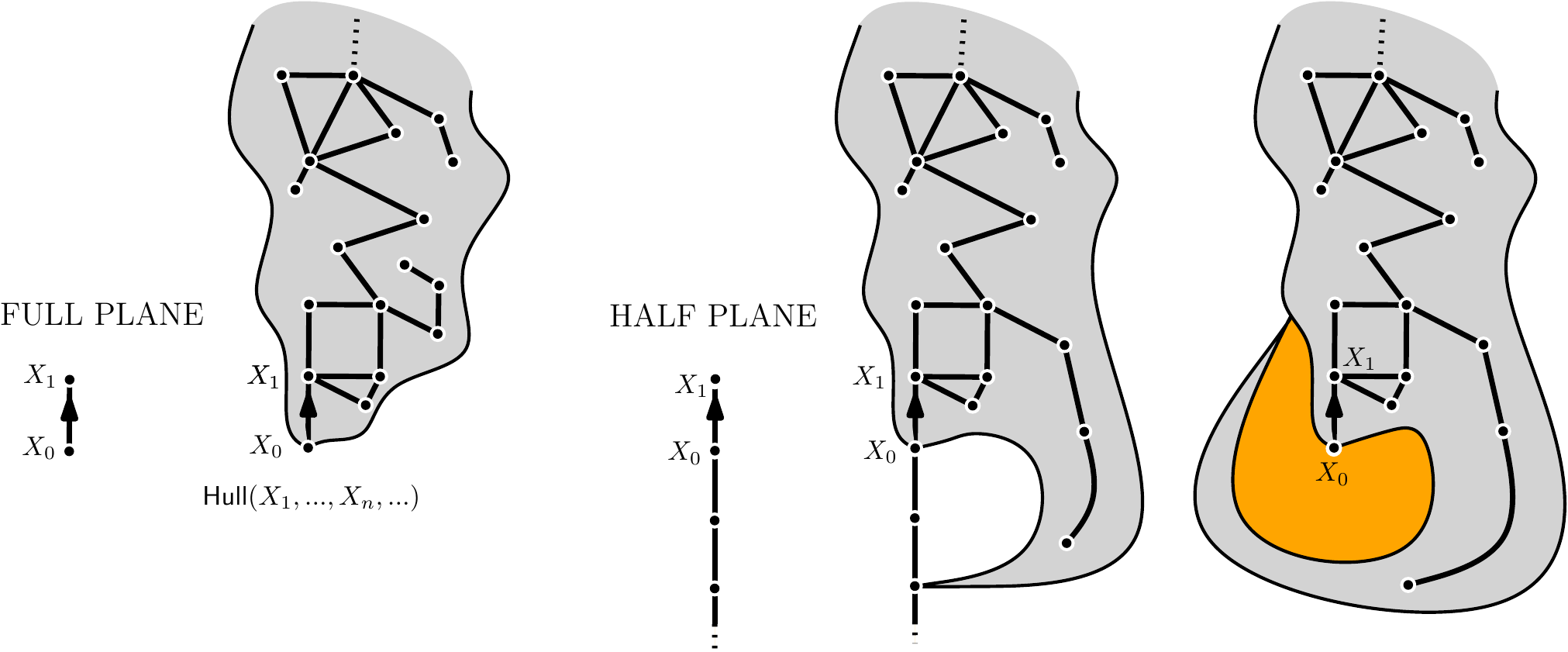}
  \caption{ \label{fig:intr2} A sequence of good instructions run in the half-plane does not intersect the infinite path except for the first two vertices.}
  \end{center}
  \end{figure}

Before proceeding to the proof of Proposition \ref{lem:intersection} let us re-interpret the last lemma in a more operative form. The history of a peeling process along a simple random walk can be seen as a sequence of random instructions, those for the walk and those for the peeling steps. The last lemma says that with positive probability, running this sequence of instructions yields a triangulation $ \mathsf{Hull}(X_{1}, ... , X_{n}, ...)$ with $X_{0}$ lying on its boundary. In this case we say that the sequence of instructions is \emph{good}. Now, imagine that we run the exact same sequence of instructions but instead of starting from the target of a single root edge, we start from the target of an infinite path whose first two vertices are $X_{1}$ and $X_{0}$, see Fig.\,\ref{fig:intr2}. In other words, we run the sequence of instruction in a half-plane. We claim that if the initial sequence of instructions is good then the hull created in the new process will not touch the infinite path except for the first two vertices.

To see this, imagine, by contradiction, that when running the sequence of instructions in the half-plane, a given peeling step reaches the infinite path further than $X_{0}$. Then, if we were in the plane, this peeling step would have gone around $X_{0}$ to reach the other side of the current explored triangulation. Doing so, it would have swallowed $X_{0}$ and thus $X_{0}$ could be on the boundary of $ \mathsf{Hull}(X_{1}, ... , X_{n}, ...)$ anymore. Contradiction.

\begin{figure}[!h]
 \begin{center}
 \includegraphics[width=13cm]{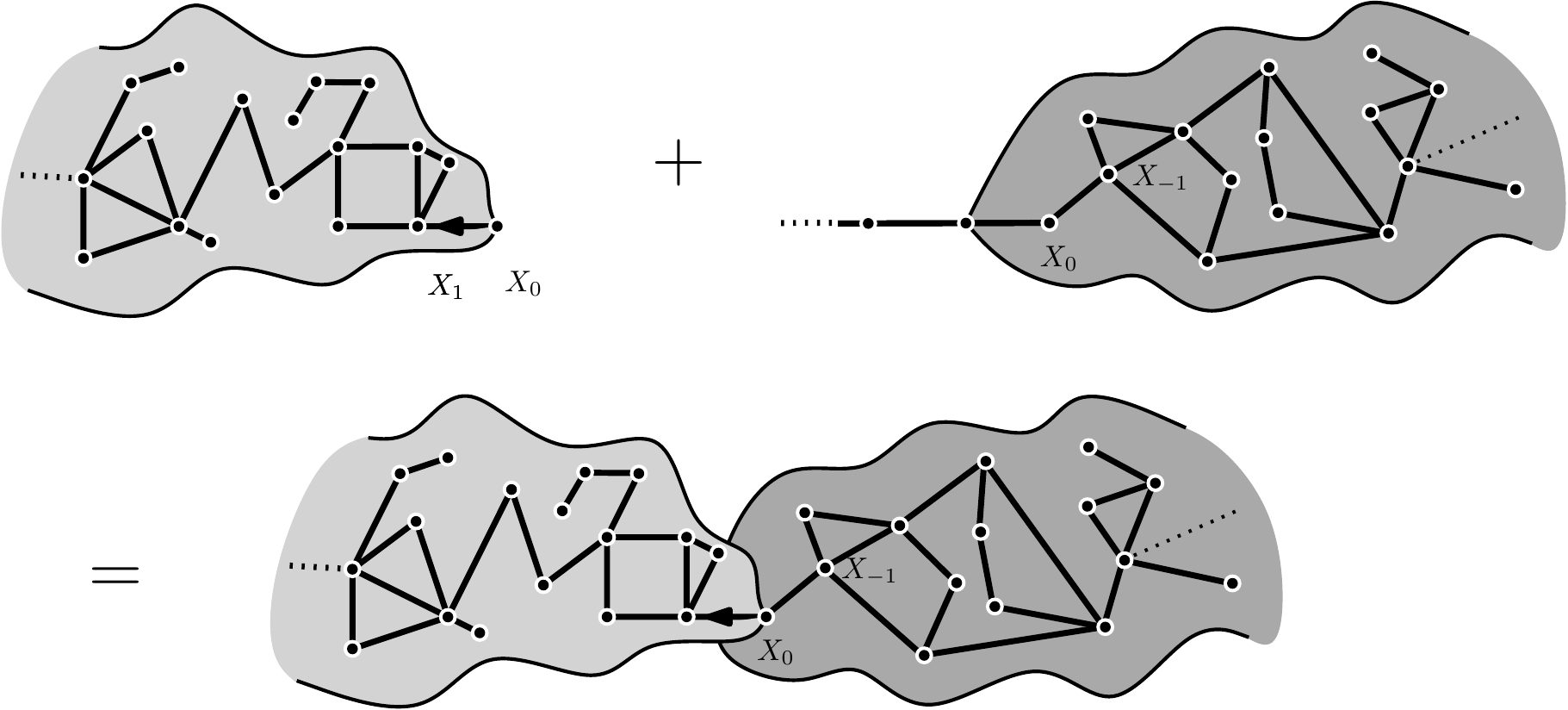}
 \caption{Illustration of the proof of Proposition \ref{lem:intersection}. \label{fig:illustr}}
 \end{center}
 \end{figure}

\proof[Proof of Proposition \ref{lem:intersection}] We give the main ideas of the proof and leave some details to the careful reader. By the last lemma, there is a positive chance that the first point $X_{0}$ lies on the boundary of $ \mathsf{Hull}(X_{1}, ... , X_{n}, ...)$, which is the triangulation explored during the peeling along the simple random walk. This implies that  $ \mathbf{T}_{\kappa} \backslash \mathsf{Hull}(X_{1}, ... , X_{n}, ... )$ is an infinite triangulation of the half-plane. A moment's though shows that the peeling process is still valid in this remaining lattice to the condition of setting $p=\infty$ in the transition probabilities. This perfectly makes sense and it can be checked (but will not be required in the argument) that this lattice has the law of Angel \& Ray's infinite triangulation of the half-plane of parameter $\alpha$ related to $\kappa$ by \eqref{eq:alphakappa}. Now, imagine that we start another independent random walk from $X_{0}$ denoted by $X_{-1},X_{-2}, ... , X_{-n},...$ and explore the rest of the lattice along it. We assume that the first step $X_{-1}$ does not belong to $\mathsf{Hull}(X_{1}, ... , X_{n}, ...)$ and so, as long as the walk does not touch $\mathsf{Hull}(X_{1}, ... , X_{n}, ...)$, this exploration can be seen as an exploration of a half-plane where $\mathsf{Hull}(X_{1}, ... , X_{n}, ...)$ has been contracted onto a half-line, see Fig.\,\ref{fig:illustr}

By the very same argument that yielded \eqref{eq:staypositive} we deduce from Lemma \ref{lem:X0} that the sequence of instructions of this new exploration is good with positive probability. That is, the walk $X_{-1}, ... , X_{-n},...$ stays  in $  \mathbf{T}_{\kappa} \backslash\mathsf{Hull}(X_{1}, ... , X_{n})$ and the only vertices in common between $ \mathsf{Hull}(X_{1}, ... , X_{n}, ...)$ and $ \mathsf{Hull}(X_{0}, X_{-1}, ... , X_{-n}, ...)$ are those within distance $1$ of $X_{0}$. In particular, on this event we have 

$$ \{X_{i}\}_{i \geq 0} \cap \{ X_{i}\}_{i \leq 0} = \{X_{0}\}.$$
By the conditional independence of the two explorations we deduce that the last event has a positive probability. Consequently, with positive probability $ \mathbf{T}_{\kappa}$ does not possess the intersection property. By ergodicity (see the proof of Lemma \ref{prop:ergodic}) almost surely $ \mathbf{T}_{\kappa}$ does not possess the intersection property. This completes the proof of Proposition \ref{lem:intersection}. \endproof

\subsection{Proof of Theorem \ref{thm:hyperbolic} \emph{via} entropy}
\label{sec:entropy}
 Combining Proposition \ref{lem:intersection} and the result of \cite{BCGLiouville} we deduce that $ \mathbf{T}_{\kappa}$ is non-Liouville a.s.\,\,when $\kappa \in (0,\frac{2}{27})$.  To finish the proof of Theorem \ref{thm:hyperbolic} is remains to prove that $s_{\kappa}>0$. For this we shall use the notion of entropy. The entropy the $n$th position of the simple random walk is the random variable defined by 
 $$ H_{n}  := \sum_{x \in \mathbf{T}_{\kappa}} \varphi \left(  \mathbb{P}(X_{n} = x)\right) \quad \mbox{ where } \quad \varphi(x) = -x \log(x).$$ Since $ \mathbf{T}_{\kappa}$ is stationary and non-Liouville  \cite[Theorem 3.2]{BCstationary}  implies that   \begin{eqnarray} \label{eq:BC} n^{-1} E[H_{n}]  &\xrightarrow[n\to\infty]{}& h > 0. \end{eqnarray} Actually, the paper \cite{BCstationary} deals with simple graphs but the proof goes through \emph{mutatis mutandis}. For technical reasons we turn this convergence in mean into an almost sure statement:
 
 \begin{lemma} We have $\limsup_{n\to \infty} n^{-1}{ H_{n}} >0$ with positive probability.
 \end{lemma}
 \proof We argue by contradiction and suppose $H_{n}/n \to 0$ almost surely. The proof of \cite[Proposition 3.1]{BCstationary} shows that $H_{n}$ is stochastically bounded by $n$ copies of (dependent) variables $H_{1,i}$ for $i \in \{1,... , n\}$ having the same law as $H_{1}$. Hence it follows by Cauchy--Schwarz inequality that 
 $$ E[H_{n}^2] \leq E\left[\left( \sum_{i=1}^n H_{1,i} \right)^2\right] = \sum_{1\leq i,j \leq n}E[H_{1,i} H_{1,j}] \leq n^2 \sqrt{E[H_{1}^2]E[H_{1}^2]}.$$
By the standard bound on the entropy we have $H_{1} \leq \log( |B_{1}( \mathbf{T}_{\kappa})|) \leq \log( |\overline{B}_{1}( \mathbf{T}_{\kappa})|)$. We leave the reader check that $|\overline{B}_{1}( \mathbf{T}_{\kappa})|$ has an exponential tail, in particular 
$$ E[H_{1}^2] \leq  E[\log^2(|\overline{B}_{1}( \mathbf{T}_{\kappa})|)] < \infty.$$ Consequently $(H_{n}/n)_{n \geq 1}$ is bounded in $ \mathbb{L}^2$ hence uniformly integrable. Since we supposed $H_{n}/n \to 0$ a.s., by dominated convergence this forces $h=0$: contradiction with \eqref{eq:BC}!
 \endproof
 
  We now adapt the proof of \cite[Proposition 3.6]{BCstationary} and demonstrate that  the last lemma  implies positive speed for the simple random walk. For this, fix $ \varepsilon>0$ and introduce the event $A^ \varepsilon_{n}= \{ \mathrm{d_{gr}}(X_{0},X_{n}) \leq ( s_{\kappa}+ \varepsilon)n\}$. To simplify notation we write $B_{r}$ for $B_{r}( \mathbf{T}_{\kappa})$. We decompose the entropy $H_{n}$ as follows \begin{eqnarray*}  \sum_{x \in \mathbf{T}_{\kappa}} \varphi \left(  \mathbb{P}(X_{n} = x)\right) &=& \sum_{x \in  B_{n( s_{\kappa}+ \varepsilon)}} \varphi(\mathbb{P}(X_{n}=x)) + \sum_{x \in B_{n} \backslash B_{n( s_{\kappa}+ \varepsilon)}} \varphi(\mathbb{P}(X_{n}=x))\\ 
   & \underset{ \varphi \mathrm{\ is \ concave}}{\leq} & \left( \sum_{x \in B_{n( s_{\kappa}+ \varepsilon)}} \mathbb{P} (X_{n}=x)\right)\log \left(  \frac{ | B_{n( s_{\kappa}+ \varepsilon)}|}{ \sum_{x \in B_{n( s_{\kappa}+ \varepsilon)}} \mathbb{P}(X_{n}=x)} \right) \\ 
     &+& \left(\sum_{x \in B_{n} \backslash B_{n( s_{\kappa}+ \varepsilon)} } \mathbb{P}(X_{n}=x) \right) \log\left( \frac{| B_{n} \backslash B_{n( s_{\kappa}+ \varepsilon)}|}{\sum_{x \in B_{n} \backslash B_{n( s_{\kappa}+ \varepsilon)} } \mathbb{P}(X_{n}=x) } \right) \\
     &= & \varphi\Big(  \mathbb{P}(A_{n}^ \varepsilon)\Big)+ \mathbb{P}(A_{n}^ \varepsilon) \log \left(|B_{n( s_{\kappa}+ \varepsilon)}| \right)  \\ &+& 
      \varphi\Big(  1-\mathbb{P}(A_{n}^ \varepsilon)\Big)+ \big(1-\mathbb{P}(A_{n}^ \varepsilon)\big) \log \left(| B_{n} \backslash B_{n(s_{\kappa}+ \varepsilon)}|\right) \end{eqnarray*}  
  We now divide by $n$ and take $\limsup_{n \to \infty}$. The left-hand side becomes positive with positive probability by the last lemma. On the other hand (the right one), from  \eqref{eq:abstractspeed} we deduce that $\mathbb{P}( A_{n}^ \varepsilon) \to 1$ almost surely under $P$ and so $\varphi( \mathbb{P}(A_{n}^ \varepsilon)) \to 0$ and $ \varphi( 1-\mathbb{P}(A_{n}^ \varepsilon)) \to 0$ as $n \to\infty$ almost surely for $P$. Also Theorem \ref{thm:volume} shows that for any $ u >0$ we have  \begin{eqnarray*} \frac{\log(|B_{un}|)}{n} & \xrightarrow[n\to\infty]{P-a.s.} & u \log\left( \frac{\alpha+ \delta_{\kappa}}{\alpha- \delta_{\kappa}}\right). \end{eqnarray*} Finally we get  
 $$ \limsup_{n \to \infty} \frac{H_{n}}{n} \leq ( s_{\kappa}+ \varepsilon)\log\left( \frac{\alpha+  \delta_{\kappa}}{\alpha- \delta_{\kappa}}\right).$$
 and conclude that $s_{\kappa}>0$ with positive probability and thus almost surely by \eqref{eq:abstractspeed}. This finishes the proof of Theorem \ref{thm:hyperbolic}.


\section{Comments} \label{sec:comments}

\subsection{Local limit of triangulations in high genus}

Recall that $ \mathcal{T}_{n,g} $ denotes the set of all (rooted) triangulations of the torus of genus $g \geq 0$ with $n$ vertices and that  $T_{n,g}$ is a random uniform element in $ \mathcal{T}_{n,g}$. Euler's formula shows that any triangulation $ t \in  \mathcal{T}_{n,g}$ has $3(n+2g-2)$ edges. Hence, when $g = [\theta n]$, the mean degree of $T_{n,[\theta n]}$ is equal to 
$$ \frac{6(n+2[\theta n]+2)}{n}  \xrightarrow[n\to\infty]{} 6(1+2 \theta).$$
However, the notion of mean degree is  \emph{not} continuous for the local topology and is not even clearly defined for an infinite triangulation. See the phenomenon appearing in the case of unicellular maps \cite[Remark 5]{ACCR13}. To get a continuous observable for the local topology, we rather look at the mean of the inverse of degree of the root vertex $\rho_{n}$ in $T_{n,g}$. Indeed, since the root vertex in $T_{n,g}$ in chosen proportionally to its degree we have 
 \begin{eqnarray*} E[ \mathrm{deg}(\rho_{n})^{-1}] &=& \frac{1}{  \# \mathcal{T}_{n,g}} \frac{1}{6(n+2g-2)} \sum_{t \in \mathcal{T}_{n,g}} \sum_{ x \in t} \mathrm{deg}(x) \cdot \frac{1}{ \mathrm{deg}(x)}\\  &=& \frac{n}{6(n+2g-2)} \quad \underbrace{\xrightarrow[n\to \infty]{} \frac{1}{6(1+ 2\theta)}}_{ \mbox{for }g=[\theta n]}.  \end{eqnarray*}
Notice that the degree of the root vertex is indeed a continuous function for the local topology. Hence, we can sharpen the conjecture stated at the end of the introduction: for $\kappa \in (0, 2/27]$, let $$f( \kappa) = E\left[ \left( \mathrm{degree \ of \ the \ origin \ in \ } \mathbf{T}_{\kappa}\right)^{-1}\right].$$
It is easy to see from the peeling construction of $ \mathbf{T}_{\kappa}$ that $f$ is continuous and satisfies $f(0^+) = 0$ and $f( \frac{2}{27}) = 1/6$ (case of the UIPT). We believe that $f$ is in fact strictly decreasing and that
\begin{conjecture}[with I. Benjamini] For any $ \theta \geq 0$, let $\kappa \in (0, \frac{2}{27}]$ be such that $ f(\kappa) = (6(1+2 \theta))^{-1}$ then we have the following convergence in distribution for the local topology
$$T_{n,[\theta n]} \quad \xrightarrow[n\to\infty]{(d)} \quad \mathbf{T}_{\kappa}.$$
\end{conjecture}
This conjecture would follow from precise enumerative formulas on $ \# \mathcal{T}_{n,g}$ when $n$ and $g$ are both tending to infinity (the known results focus on asymptotics as $n \to \infty$ and then $g \to \infty$, \cite{GJ08}), see the arguments in \cite{ACCR13}.
\subsection{Perspectives} First of all, let us mention that we restricted ourselves to $2$-connected triangulations mainly to take advantage of the calculations already performed by Angel \& Ray in \cite{AR13} and by Ray in \cite{Ray13}. This whole work could be extended to $1$-connected triangulations or other types of maps (e.g.\,\,quadrangulations) to the price of adapting the constants.  

Also, it is likely that site percolation on $ \mathbf{T}_{\kappa}$ can be treated by similar means as in \cite{Ray13} and would yield almost identical results. It is pretty clear that in the hyperbolic regime $ \mathbf{T}_{\kappa}$ do not admit any scaling limits in the Gromov--Hausdorff sense (the hull $ \overline{B}_{r}( \mathbf{T}_{\kappa})$ contains a number of tentacles reaching distance $2r$ that tends to infinity as $r \to \infty$), however its conformal structure might be of interest. Finally, the geometric relations (underlying the proof of Proposition \ref{lem:intersection}) between the half-planar lattices of \cite{AR13} and those defined in this work deserve to be explored in more details.  \bigskip

\bibliographystyle{siam}

\end{document}